\documentclass[12pt]{amsart}
\usepackage{amsmath,amssymb,amsfonts,amsthm}
\usepackage[hidelinks]{hyperref}
\hypersetup{colorlinks=true,linkcolor=blue, citecolor=blue, urlcolor=blue}
\usepackage{tikz}
\usepackage{tikz-cd}
\usepackage{enumerate}
\usepackage{graphicx, array}
\usepackage{float}
\usepackage{xfrac} 

\usepackage{graphicx}
\graphicspath{{figures/}}
\usepackage[labelfont=small,labelformat=simple]{subcaption}

\theoremstyle{plain}

\theoremstyle{definition}
\newtheorem{definition}{Definition}[section]

\theoremstyle{remark}
\newtheorem{remark}{Remark}[section]

\usepackage[T5,T1]{fontenc}

\usepackage[
backend=biber,
style=alphabetic,
sorting=nyt
]{biblatex}
\DeclareFieldFormat*{title}{\textit{#1}}
\DeclareFieldFormat{journaltitle}{\textnormal{#1}}
\DeclareListFormat{publisher}{\textnormal{#1}}

\usepackage[left=2.5cm, right=2.5cm,top=3cm,bottom=3cm]{geometry}

\newcommand{\Z}{\mathbb{Z}}
\newcommand{\R}{\mathbb{R}}
\newcommand{\C}{\mathbb{C}}

\newcommand{\abs}[1]{\left\vert{#1}\right\vert}

\newcommand{\del}[2][]{\frac{\partial{#1}}{\partial{#2}}}

\newcommand{\vungoc}{V{\~u} Ng{\d{o}}c}
\title{Towards hypersemitoric systems}

\author{Tobias V\r{a}ge Henriksen}
\address{Tobias V\r{a}ge Henriksen: 
  Bernoulli Institute for Mathematics\\ 
  Computer Science and Artificial Intelligence\\
  University of Groningen\\ 
  P.O. Box 407, 9700 AK Groningen, The Netherlands
}
\email{t.v.henriksen@rug.nl}

\author{Sonja Hohloch}
\address{Sonja Hohloch: 
  Department of Mathematics\\ 
  University of Antwerp\\ 
  Middelheimlaan 1, B-2020 Antwerp, Belgium
}
\email{sonja.hohloch@uantwerpen.be}

\author{Nikolay N. Martynchuk}
\address{Nikolay N. Martyncuk: 
  Bernoulli Institute for Mathematics\\ 
  Computer Science and Artificial Intelligence\\
  University of Groningen\\ 
  P.O. Box 407, 9700 AK Groningen, The Netherlands
}
\email{n.martynchuk@rug.nl}

\begin{document}

\date{\today}
\maketitle

\begin{abstract}
This survey gives a short and comprehensive introduction to a class of finite-dimensional integrable systems known as \textit{hypersemitoric systems}, recently introduced by Hohloch and Palmer in connection with the solution of the problem how to extend Hamiltonian circle actions on symplectic 4-manifolds to integrable systems with `nice' singularities. 
The quadratic spherical pendulum, the Euler and Lagrange tops (for generic values of the Casimirs), coupled-angular momenta, and the coupled spin oscillator system are all examples of hypersemitoric systems.
Hypersemitoric systems are a natural generalization of so-called semitoric systems (introduced by V{\~u} Ng{\d{o}}c) which in turn generalize toric systems. Speaking in terms of bifurcations, semitoric systems are `toric systems with/after supercritical Hamiltonian-Hopf bifurcations'. Hypersemitoric systems are `semitoric systems with, among others, subcritical Hamiltonian-Hopf bifurcations'.
Whereas the symplectic geometry and spectral theory of toric and semitoric sytems is by now very well developed, the theory of hypersemitoric systems is still forming its shape. This short survey introduces the reader to this developing theory by presenting the necessary notions and results as well as its connections to other areas of mathematics and mathematical physics.
\end{abstract}
\section{Introduction}

Integrable Hamiltonian systems play an important role in mathematical and physical sciences. For instance, within celestial mechanics, there is the Kepler problem, and, within quantum mechanics, there is the Jaynes-Cummings model, which are both integrable. Integrable systems are very special dynamical systems exhibiting regular (as opposed to chaotic) behaviour in the sense that there exist a maximal number of (independent, see Definition \ref{def:IntegrableSystem}) integrals of motion, allowing one to at least in principle integrate the equations of motion.

Dynamics of a finite-dimensional integrable Hamiltonian system, defined by means of a proper momentum map (see Definition \ref{def:IntegrableSystem}), is generically constrained to $n$-dimensional tori, where $n$ is the number of degrees of freedom. These tori turn out to be Lagrangian submanifolds of the underlying symplectic manifold on which the Hamiltonian system is defined, and thus an integrable system can be seen as a singular Lagrangian torus fibration over a certain subset of $\R^n$, see in particular the papers by Mineur \cite{Mineur1936}, Arnol'd \cite{Arnold1963}, Weinstein \cite{Weinstein1971} and Duistermaat \cite{Duistermaat1980}. This motivates one to study integrable systems using techniques from symplectic geometry. 

The singular fibres of these singular Lagrangian torus fibrations reflect a non-trivial geometric or dynamical property of the underlying integrable system. The most prominent examples being the monodromy around a focus-focus point and bifurcations of Liouville tori, which we will address below. 

In the context of symplectic classification of integrable systems it is known how to classify a number of different types of such (`typical') singularities: a saddle singularity (in one degree of freedom) by Dufour, Molino, and Toulet \cite{Dufour1994}, an elliptic singularity (in any dimension) by Eliasson \cite{Eliasson1984, Eliasson1990}, a focus-focus singularity (in dimension 2) by \vungoc\ \cite{Vu-Ngoc2003}, and a parabolic singularity by Bolsinov, Guglielmi, and Kudryavtseva \cite{Bolsinov2018} and Kudryavtseva and Martynchuk \cite{Kudryavtseva2021}. See also the recent breakthrough results concerning symplectic classification in the real-analytic category by Kudryavtseva \cite{Kudryavtseva2020} and by Kudryavtseva and Oshemkov \cite{Kudryavtseva2022}. 

In the context of global classification of integrable systems, Pelayo and V\~u Ng{\d o}c \cite{Pelayo2009} showed that a large class of physically important systems known as semitoric systems are classified by a set of 5 invariants. This is one of the few known explicit results in the global symplectic classification of integrable systems, apart from the classical Delzant's \cite{Delzant1988} construction and the work of Zung \cite{Zung2003} relating the semi-local (i.e.\ in a neighbourhood of a singular fibre) and global classification problems. We refer to Sections~\ref{sec:semitoric_overview} and \ref{sec:semitoric} for more details on semitoric systems. 

What is currently missing in the literature is a detailed discussion of systems beyond semitoric type: Whereas the \textit{topological}
classification of such systems is a well developed theory going back to Duistermaat and Fomenko and Zieschang (see e.g.\ Bolsinov and Fomenko \cite{Bolsinov2004} and the references therein), a more refined (e.g.\ symplectic) analysis is currently an open problem for in fact the majority of such systems. In particular, what is missing is a detailed analysis of a generalisation of semitoric systems additionally allowing hyperbolic-regular,  hyperbolic-elliptic, and parabolic points,  known as \textit{hypersemitoric} systems. The latter class was introduced by Hohloch and Palmer 
\cite{Hohloch2021} in connection with the problem of extending Hamiltonian circle actions on symplectic 4-manifolds to integrable systems, which they solved within this class of systems, see Hohloch and Palmer \cite{Hohloch2021} for details.
Hypersemitoric systems thus present a challenging platform for the further study by both geometers and analysists and this survey is devised as a quick introduction. 

Nevertheless, note that the class of hypersemitoric systems does not include all possible singularities that may arise in 4-dimensional integrable systems: the underlying global $S^1$-action prevents the existence of hyperbolic-hyperbolic singularities; moreover, the definition of hypersemitoric systems excludes most of the `typical' degenerate $S^1$-invariant singularities, see Kalashnikov's \cite{Kalashnikov1998} list. There exists another class of integrable systems, namely hyperbolic semitoric systems (cf.\ \cite[Definition 3.2]{DullinPelayo2016}), which, if one considers the union with semitoric systems, contains hypersemitoric system, see Remark \ref{rem:hyperbolic_semitoric}. The hyperbolic semitoric systems do \emph{include} all `typical' degenerate $S^1$-invariant singularities in Kalashnikov's \cite{Kalashnikov1998} list.

\subsection*{Organization of the paper}

The rest of this paper is organized as follows: In Section \ref{sec:overview}, we give the definition of (Liouville) integrability, before defining toric, semitoric, and hypersemitoric systems. Moreover, we explain some important properties of integrable systems and give a short survey over the theory of atoms and molecules. In Section \ref{sec:semitoric}, we discuss semitoric systems in detail, i.e., their symplectic classification in terms of five invariants and how one may obtain a semitoric system from a toric one. Eventually, we recall some important examples. In Section \ref{sec:hypersemitoric}, we consider hypersemitoric systems: we first discuss flaps and pleats, which occur in the momentum image of hypersemitoric systems. Then we consider how one may obtain hypersemitoric systems from (semi)toric systems before we briefly explain an explicit example.

\subsection*{Acknowledgements}

The authors are very grateful to \'Alvaro Pelayo and San V{\~u} Ng{\d{o}}c for useful comments and suggestions that helped to improve the original version of this work.
The first author was fully supported by the {\em Double Doctorate Funding} of the Faculty of Science and Engineering of the University of Groningen.
Moreover, all authors were partially supported by the \textit{FNRS-FWO Excellence of Science (EoS) project `Symplectic Techniques in Differential Geometry' G0H4518N}. 

\section{Definitions, conventions, and background}\label{sec:overview}

In this section, we give an outline of integrability with an emphasis on integrable systems defined on 4-manifolds and admitting a global effective Hamiltonian circle action. Hypersemitoric systems are a certain class of systems of this type.  We start by recalling the classical Arnol'd-Liouville-Mineur theorem, and then move from toric to semitoric to hypersemitoric systems. We also show how the theory relates to the general frameworks of monodromy and bifurcations of Liouville tori, i.e., Fomenko-Zieschang theory. 

\subsection{Integrable systems}
\label{sec:IntSys}

Let $(M, \omega)$ be a \emph{symplectic manifold} of dimension $2n$. Since the symplectic form is non-degenerate, for any function $f \in C^{\infty}(M,\R)$, there exists a unique vector field $X_{f}$, called the \emph{Hamiltonian vector field} of $f$, such that $\iota_{X_{f}} \omega = - df$. The function $f$ is called the \emph{Hamiltonian}, and $\dot{z} = X_{f}(z)$ is called a \emph{Hamiltonian system}, sometimes briefly denoted by $X_{f}$. For two Hamiltonians $f,g \in C^{\infty}(M,\R)$, the \emph{Poisson bracket} is defined by $\{f, g\} := \omega(X_{f}, X_{g})$. If $\{f, g\} = 0$, then $f$ and $g$ are said to \emph{Poisson commute}. Note that $\{f, g\} = X_{f}(g)$. If $f$ and $g$ Poisson commute, then $g$ is called a \emph{(first) integral} of $X_{f}$.

\begin{definition}\label{def:IntegrableSystem}
A Hamiltonian system $X_{H}$ on a $2n$-dimensional symplectic manifold $(M, \omega)$ is said to be \textit{completely integrable} (or briefly \emph{integrable}) if there exist $n$ functionally independent integrals $f_{1} := H, f_{2},\dots,f_{n}$ of $X_{H}$, i.e.\ their gradients are almost everywhere linearly independent on $M$, the integrals all Poisson commute with each other, and the flows of $X_{f_1}$, \dots, $X_{f_n}$ are complete.
 A shorter notation is $(M, \omega, F=(f_{1},\dots,f_{n}))$ and $F$ is often referred to as the \emph{momentum} or \emph{integral} map of the system.
\end{definition}

A point $p\in M$ is {\em regular} if the rank of $DF_p$ is maximal and {\em singular} otherwise.  A value of $F$ is {\em regular} if all points in the preimage are regular, and {\em singular} otherwise.
Similarly, one defines what it means for a fibre $F^{-1}(r)$ of $F$ to be regular, resp., singular and for a  {\em leaf} of $F$, i.e. a connected component of a fibre, to be regular, resp. singular.

The Arnol'd-Liouville-Mineur theorem \cite{Arnold1978,Mineur1936} describes the regular leaves of the foliation generated by the momentum map of a $2n$-dimensional integrable system. Each regular leaf is a Lagrangian submanifold, and if the leaf is connected and compact, then it is diffeomorphic to an $n$-torus $T^{n}$. Such a foliation will be called a \emph{Lagrangian torus fibration}. Let $r \in \R^n$ be a regular value for the momentum mapping $F$, and let $F^{-1}(r)$ be a connected and compact fibre, and hence diffeomorphic to $T^{n}$, and let $U$ be a tubular neighbourhood of $F^{-1}(r)$. The Arnol'd-Liouville-Mineur theorem also tells us that $U$ is diffeomorphic to $V \times T^{n}$, where $V$ is an open set of $\R^{n}$. On $V \times T^n$, there exists coordinates $I_{1}, \dots, I_{n}, \phi_{1}, \dots, \phi_{n}$, called \emph{action-angle coordinates}. Here each $I_i$ for $i = 1, \dots, n$ is a function of the $f_{i}$'s, whilst each $\phi_{i}$ is a standard angle coordinate on $T^{n}$. In action-angle coordinates, the symplectic form becomes $\omega = \sum d\phi_{i} \wedge dI_{i}$. Note that, in general, action-angle coordinates only exist locally. Duistermaat \cite{Duistermaat1980} showed that there can exist obstructions to the global existence of action-angle coordinates in terms of the (Hamiltonian) monodromy and the Chern class on the topological level as well as the Lagrangian class on the symplectic level. 

For us, monodromy will play an essential role so that we will recall its definition here; for more detail see \cite{Duistermaat1980}. Let $F : M \to B$ be a Lagrangian torus fibration over an $n$-dimensional manifold $B$ and denote by $R \subseteq B$ the set of the regular values of $F$. 
Then there exists a natural covering
$$\bigcup_{r \in R} \textup{H}_1(F^{-1}(r)) \to R,$$
where $\textup{H}_1(F^{-1}(r))$ is the first homology group of $F^{-1}(r)$ with integer coefficients.
Because of this, there is a natural representation of $\pi_{1}(R)$ into the group $\textnormal{SL}(n, \Z)$ of automorphisms of the lattice 
$\textup{H}_1(F^{-1}(r)) \simeq \mathbb Z^n.$ This representation is called the \emph{Hamiltonian monodromy} of $F : M \to B$ (or of $F : M \to R$). Thus, to any loop $\gamma$ in $R$, one can assign an $n\times n$ integer matrix called the \emph{monodromy} or the \emph{monodromy matrix} along $\gamma$.

Note that  Lagrangian torus fibrations are allowed to have singular points and these are precisely the points that encode essential properties of the underlying integrable system. One has in particular been interested in \emph{non-degenerate singular points}, i.e.\ points for which the Hessians of the integrals span a Cartan subalgebra in the real symplectic Lie algebra $\textnormal{sp}(2n, \R)$ (cf.\ Bolsinov and Fomenko \cite{Bolsinov2004}). Locally one can describe such singularities by local normal forms (cf., among other, the works by Eliasson \cite{Eliasson1984, Eliasson1990}, Miranda and Zung \cite{Miranda2004}, and \vungoc\ and Wacheux \cite{Vu-Ngoc2013}): in a neighbourhood $U$ of a non-degenerate singular point, one can find local symplectic coordinates $(x_{1}, \dots, x_{n}, \xi_{1}, \dots, \xi_{n})$ such that the symplectic form takes the form $\omega = \sum_{i=1}^{n} dx_{i} \wedge d\xi_{i}$ in $U$, and $n$ functionally independent smooth integrals $q_{1}, \dots, q_{n} : U \to \R$ Poisson commuting with all  $f_{1}, \dots, f_n$ such that $q_{i}$ is one of the following possible components:
\begin{itemize}
    \item regular component: $q_{i} = x_{i}$,
    \item elliptic component: $q_{i} = \frac{1}{2}(x_{i}^{2} + \xi_{i}^{2})$,
    \item hyperbolic component: $q_{i} = x_{i}\xi_{i}$,
    \item focus-focus components (exist in pairs): $q_{i} = x_{i}\xi_{i} + x_{i+1}\xi_{i+1}$ and $q_{i+1} = x_{i}\xi_{i+1} - x_{i+1}\xi_{i}$.
\end{itemize}
We will eventually focus on $4$-dimensional integrable systems. In that case, the following six different types of non-degenerate singular points can occur:
\begin{itemize}
    \item rank 0: elliptic-elliptic, hyberbolic-hyperbolic, elliptic-hyperbolic and focus-focus,
    \item rank 1: elliptic-regular and hyperbolic-regular.
\end{itemize}
Williamson \cite{Williamson1936} (see also Bolsinov and Fomenko \cite[Section 1.8]{Bolsinov2004}) showed that to determine the type of a non-degenerate rank 0 singular point of a $4$-dimensional integrable system $(M, \omega, F=(f_1, f_2))$, it is sufficient to find the eigenvalues for the Hessian of the linear combination $c_1 f_1 +  c_2 f_2$ for generic $c_1, c_2 \in \R$ at this singular point since
\begin{itemize}
    \item elliptic components have pairs of purely imaginary eigenvalues,
    \item hyperbolic components have pairs of purely real eigenvalues,
    \item focus-focus components have quadruples of complex eigenvalues with non-zero real- and imaginary parts.
\end{itemize}
Note also that, if $\lambda$ is an eigenvalue of multiplicity $k$, then so are $-\lambda$, $\Bar{\lambda}$, and $-\Bar{\lambda}$ (cf.\ van der Meer \cite[Proposition 1.28]{vanderMeer1985}).

Concerning monodromy, we note that if $\Lambda$ 
is a (compact) leaf containing $n$ singular points of which all are of focus-focus type, then it has been shown that the monodromy around $\Lambda$ is given by 
\begin{equation*}
    M = \begin{pmatrix} 1 & n \\ 0 & 1 \end{pmatrix},
\end{equation*}
 see the works by Matsumoto \cite{Matsumoto1989}, Lerman and Umanskii \cite{Lerman1994}, Matveev \cite{Matveev1996}, and Zung \cite{Zung1997}.
This result will be drawn on again in our discussion of semitoric and hypersemitoric systems.

\subsection{Toric systems}

Let us start with the `easiest' class of integrable systems:

\begin{definition}
Let $(M,\omega,F)$ be an integrable system with $M$ compact and connected. If all integrals of $(M,\omega,F)$ generate an effective $S^{1}$-action, then the system is said to be a \emph{toric system}. 
\end{definition}

Atiyah \cite{Atiyah1982} and Guillemin and Sternberg \cite{Guillemin1982} showed that the image of the momentum map of a toric system is a convex polytope, called the \emph{momentum polytope}. Later, Delzant \cite{Delzant1988} showed that toric systems are classified up to isomorphism by their momentum polytope. Delzant's classification was then extended to non-compact manifolds by Karshon and Lerman \cite{Karshon2015}. Note that the singular points of a toric system are all non-degenerate and only contain components of elliptic or regular type.

\subsection{Semitoric systems} \label{sec:semitoric_overview}

Delzant's \cite{Delzant1988} classification of toric manifolds has been generalized by Pelayo and V{\~u} Ng{\d{o}}c \cite{Pelayo2009, Pelayo2011} together with Palmer and Pelayo and Tang \cite{Palmer2019} to the following class of integrable systems, called ``semitoric systems''. Semitoric systems are a natural class of systems, generalizing toric systems by relaxing the assumption of periodicity on one of the integrals defining the system. Semitoric systems are closely related to so called \emph{almost-toric system}, see for instance Symington \cite{Symington2002} and V{\~u} Ng{\d{o}}c \cite{Vu-Ngoc2007}. The notion ``semitoric'' is natural, and has been used in different contexts, including symplectic geometry of Hamiltonian torus action by Karshon and Tolman \cite{Karshon2001}, integrable systems V{\~u} Ng{\d{o}}c \cite{Vu-Ngoc2007} and Pelayo and V{\~u} Ng{\d{o}}c \cite{Pelayo2009, Pelayo2011}, partially equivariant embedding problems in toric geometry by Pelayo \cite{Pelayo2007}, and mathematical physics by Martini and Taylor \cite{Martini2015}. We refer to Pelayo \cite[Sections 8.2, 8.3 and 9.3]{pelayo2023} for further discussion and references.

\begin{definition}[Pelayo and V{\~u} Ng{\d{o}}c {\cite[Definition 2.1]{Pelayo2009}}]\label{def:SemitoricSystem}
Let $(M, \omega, F=(J,H))$ be a $4$-dimensional integrable system, where $M$ is connected. Then $(M, \omega, F=(J,H))$ is a \emph{semitoric system} if
\begin{enumerate}
    \item $J$ is proper and generates an effective $S^{1}$-action,
    \item $F$ has only non-degenerate singularities (if any) and none of them admit hyperbolic components.
\end{enumerate}
\end{definition}

Note that, under the assumptions of Definition \ref{def:SemitoricSystem}, V{\~u} Ng{\d{o}}c \cite[Theorem 3.4]{Vu-Ngoc2007} showed that the fibres of $F$ are connected, thus generalizing the connectivity statement from the toric case as shown by Atiyah \cite{Atiyah1982} and Guillemin and Sternberg \cite{Guillemin1982}.

The main difference between toric and semitoric systems is the possible appearance of focus-focus singular points. Note that if $c \in F(M)$ is a focus-focus singular value, then its preimage $F^{-1}(c)$ has the shape of a so-called \emph{pinched torus} where the number of pinches equals the number of focus-focus points in the fibre, cf.\ for instance Bolsinov and Fomenko \cite{Bolsinov2004}.

V{\~u} Ng{\d{o}}c \cite{Vu-Ngoc2007} showed that one can associate an equivalence class of polygons with the image of the momentum map of a semitoric system. But unlike to the toric case, this is not enough to classify semitoric systems. Pelayo and V{\~u} Ng{\d{o}}c \cite{Pelayo2009, Pelayo2011} were able to classify so-called {\em simple} semitoric systems, i.e.\ semitoric systems for which each fibre of $J$ contains at most one focus-focus point, by formulating the following five invariants:
\begin{enumerate}[(i)]
    \item the number of focus-focus points,
    \item the Taylor series or singularity type invariant,
    \item the polygon invariant,
    \item the height invariant, and
    \item the twisting index invariant.
\end{enumerate}
Palmer, Pelayo and Tang \cite{Palmer2019} extended the result to the non-simple case, building on the symplectic classification of multi-pinched focus-focus fibres by Pelayo and Tang \cite{Pelayo2022}.

The five invariants will be discussed further in Section \ref{sec:semitoric}, where also two examples will be covered, namely the coupled angular momenta (Section \ref{sec:semitoric_cam}), and an example for which the polygon takes the shape of an octagon (Section \ref{sec:semitoric_octagon}). Other important examples of semitoric systems are the spherical pendulum (cf.\ Dullin \cite{DullinSpherPendulum}) and the Jaynes-Cummings model (cf.\ Babelon, Cantini and Dou{\c{c}}ot \cite{Babelon2009}, Pelayo and V{\~u} Ng{\d{o}}c \cite{Pelayo2012}, and Alonso, Dullin and Hohloch \cite{Alonso2019a}).

\subsection{Hypersemitoric systems}

Hohloch and Palmer \cite{Hohloch2021} considered a yet more general class of integrable systems than semitoric systems by allowing for singular points with hyperbolic components and certain degenerate singular points, namely so-called parabolic singular points: a singular point $p$ of an integrable system $(M, \omega, F=(f_{1},f_{2}))$ is {\em parabolic} if there exists a neighbourhood $U \subset M$ of $p$ with (generally non-canonical) coordinates $(x, y, \lambda, \phi)$ and functions $q_{i} = q_i(f_1,f_2)$ for $ i \in \{ 1,2\}$ of the form
$$q_{1} = x^{2} - y^{3} + \lambda y \ \mbox{ and } \  q_{2} = \lambda.$$ 
 A coordinate free definition is given in Bolsinov, Guglielmi and Kudryavtseva \cite{Bolsinov2018}. 
 Note that the same normal form in fact applies to parabolic \textit{orbits}, which means that from the smooth point of view, there is only one type of degenerate
singularities appearing in hypersemitoric systems (for more details, see Kudryavtseva and Martynchuk \cite[Theorem 3.1]{Kudryavtseva2021b}).
Parabolic points are also known under the name of \emph{cusps} or \emph{cuspidal points}. Moreover, parabolic points naturally appear as transition points between (families of) elliptic-regular and hyperbolic-regular points.

The following definition generalizes the natural notions of toric and semitoric systems we have seen earlier in this paper, and appears in recent work by Hohloch and Palmer \cite{Hohloch2021}, following also work by Kalashnikov \cite{Kalashnikov1998} as explained below.

\begin{definition}[Hohloch and Palmer {\cite[Definition 1.6]{Hohloch2021}}]\label{def:HypersemitoricSystem}
A $4$-dimensional integrable system $(M, \omega, F=(J,H))$ is called \emph{hypersemitoric} if 
\begin{enumerate}
    \item $J$ is proper and generates an effective $S^{1}$-action,
    \item all degenerate singular points of $F$ (if any) are of parabolic type.
\end{enumerate}
\end{definition}

Note that the existence of a global $S^{1}$-action prevents the appearance of hyperbolic-hyperbolic singularities in a hypersemitoric system. The original motivation for introducing this class, however, comes from the result of Hohloch and Palmer \cite[Theorem 1.7]{Hohloch2021} stating that any $4$-dimensional Hamiltonian system $X_{J}$ which generates an effective $S^1$-action is extendable to a hypersemitoric system $(M, \omega, (J,H))$. Furthermore, the set of hypersemitoric systems is open in the set of $4$-dimensional integrable systems with a global effective Hamiltonian circle action (see Kalashnikov \cite{Kalashnikov1998}).

\begin{remark}\label{rem:hyperbolic_semitoric}
Dullin and Pelayo \cite{DullinPelayo2016} showed that, starting with a semitoric system, one can use a subcritical Hamiltonian-Hopf bifurcation (which transforms a focus-focus point to an elliptic-elliptic point, see Sections \ref{sec:supercritical_Hamiltonian-Hopf} and \ref{sec:subcritical_Hamiltonian-Hopf}) to generate a flap (see Section \ref{sec:flaps_pleats}) on said system, thus creating a \emph{hyperbolic semitoric system} (cf.\ \cite[Definition 3.2]{DullinPelayo2016}). Although the name of this type of system is very similar to the name hypersemitoric, they are defined differently. Hyperbolic semitoric systems requires the same conditions as hypersemitoric systems for the integral $J$ generating a circle action. However, the set of hyperbolic singularities in hyperbolic semitoric systems are required to be non-empty, and the set of degenerate singularities is required to be isolated, not necessarily of parabolic type. Nevertheless, many hypersemitoric systems can thus be generated by performing subcritical Hamiltonian-Hopf bifurcations, together with so-called \emph{blow-ups} (also known as corner chops, see for instance Holoch and Palmer \cite{Hohloch2021} and references therein) on the (newly generated) elliptic-elliptic points.
\end{remark}

\subsection{Topological invariants: atoms and molecules}
\label{sec:AtomsMolecules}

Finally, we will recall a complete topological invariant for a generic isoenergy level of a two degree of freedom integrable system which was introduced by Fomenko and Zieschang \cite{Fomenko1990}. This invariant is intimately linked to hyperbolic-regular and elliptic-regular points
and naturally appears in (hyper)semitoric systems as well as in systems without a global $S^1$-action, which in fact form a majority of known integrable systems (including the Kovalevskaya top and many other integrable cases in rigid body dynamics, various geodesic flows, billiards, etc.).  We will follow the presentation of Bolsinov and Fomenko \cite[Sections 2, 3 and 4]{Bolsinov2004}.

Let $f$ be a Morse function on a manifold $M$. Note that the leaves of $f$ foliate the manifold. Let $x \sim y$ if and only if $x$ and $y$ are in the same leave of $f$ and denote by $\Gamma := M / \sim$ the {\em space of leaves} of $f$. 
Since $f$ is a Morse function $\Gamma$ is in fact a graph, called the \emph{Reeb graph} of $f$ on $M$ where singular leaves give rise to the vertices. There are two types of vertices:
\begin{enumerate}
    \item a vertex is called an \emph{end vertex} if it is the end of one edge only,
    \item otherwise it is called an \emph{interior vertex}.
\end{enumerate}
Note that the end vertices of a Reeb graph correspond to local minima and maxima (thus elliptic points) of the Morse function, whilst the interior vertices correspond to saddle-points (thus hyperbolic points).

Let $f \colon M \to \mathbb R$ be a Morse function on a 2-dimensional surface $M$. An \emph{atom} is a tubular neighbourhood denoted by $P^{2}$ of a singular fibre $f^{-1}(c)$ together with the fibration $f \colon P^2 \to \mathbb R$ on this neighbourhood. The atom is \emph{orientable} if the surface $P^{2}$ is orientable and \emph{non-orientable} otherwise. We now give a brief overview of the so-called \emph{simple} atoms, which are atoms whose singular fibres contain only one singular point and which are referred to as \emph{atom $A$}, \emph{atom $B$} and \emph{atom $\Tilde{B}$}. There exist many more atoms, which are defined similarly to the aforementioned ones. A more detailed exposition can be found in Bolsinov and Fomenko \cite[Section 2.4]{Bolsinov2004}.

\begin{figure}[ht]
  \centering
  \begin{subfigure}[b]{0.35\linewidth}
  \centering
  \begin{tikzpicture}[baseline]
\draw[->, line width=0.5mm] (3,0.8) node[anchor=north]{$A$} -- (3,1.3);
\draw[-, line width=0.5mm] (3,1.3) -- (3,1.8);

\draw[->, line width=0.5mm] (5,0.8) -- (5,1.3);
\draw[-, line width=0.5mm] (5,1.3) -- (5,1.8) node[anchor=south]{$A$};
\end{tikzpicture}

\caption{ }
\end{subfigure}
\hspace{1.8cm}
  \begin{subfigure}[b]{0.45\linewidth}
  \centering
  \begin{tikzpicture}
\draw[->, line width=0.5mm] (0,0) -- (0,0.5);
\draw[-, line width=0.5mm] (0,0.5) -- (0,1);
\draw[->, line width=0.5mm] (0,1.6) node[anchor=north]{$B$} -- (-0.25,2.1);
\draw[-, line width=0.5mm] (-0.25,2.1) -- (-0.5,2.6);
\draw[->, line width=0.5mm] (0,1.6) node[anchor=north]{$B$} -- (0.25,2.1);
\draw[-, line width=0.5mm] (0.25,2.1) -- (0.5,2.6);

\draw[->, line width=0.5mm] (1.5,0) -- (1.75,0.5);
\draw[-, line width=0.5mm] (1.75,0.5) -- (2,1);
\draw[->, line width=0.5mm] (2.5,0) -- (2.25,0.5);
\draw[-, line width=0.5mm] (2.25,0.5) -- (2,1);
\draw[->, line width=0.5mm] (2,1.6) node[anchor=north]{$B$} -- (2,2.1);
\draw[-, line width=0.5mm] (2,2.1) -- (2,2.6);

\draw[->, line width=0.5mm] (4,0) -- (4,0.5);
\draw[-, line width=0.5mm] (4,0.5) -- (4,1);
\draw[->, line width=0.5mm] (4,1.6) node[anchor=north]{$\Tilde{B}$} -- (4,2.1);
\draw[-, line width=0.5mm] (4,2.1) -- (4,2.6);
\end{tikzpicture}
\caption{ }
\end{subfigure}
\caption{Subfigure (a): The left graph displays a minimum, while the right graph displays a maximum. Subfigure (b): The graph to the left displays, for example, the level sets of the height function on a torus near the bottom where the level sets change from one to two circles when passing through the lower saddle. Similarly the graph in the middle displays the case at the upper saddle of the torus where the level sets change from two to one circle. These are the atoms of type $B$. The graph to the right displays a non-orientable saddle, and it is an atom of type $\Tilde{B}$.}
\label{fig:atoms}
\end{figure}
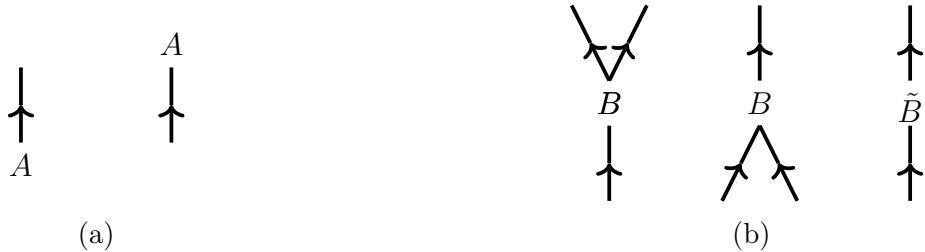

Let us first consider atom $A$, which represents the case of local minima or maxima of the function $f$. The Reeb graph of the atom is a line segment illustrating the energy levels of $f$ together with an arrow pointing in the direction of increasing energy, and a symbol $A$ illustrating the extrema. Thus, there exist two atoms of type $A$ of which the associated Reeb graphs are sketched in Figure \ref{fig:atoms}.

\renewcommand\thesubfigure{(\alph{subfigure})}

One can do a similar construction for saddles. Note, however, that there exist both orientable and non-orientable saddles, and they lead to atoms of type $B$ and $\Tilde{B}$, respectively. One can generate such atoms by considering a cylinder and gluing a strip to one of its ends (more specifically, attaching an index-1 handle). If the strip is not twisted, this can be deformed to an orientable saddle, whilst if it is twisted, it can be deformed to a non-orientable saddle. Figure \ref{fig:atoms} shows the Reeb graphs of these atoms.

There also exist atoms with more than one singular point in the singular fibre (cf.\ Bolsinov and Fomenko \cite[Section 2.6]{Bolsinov2004}). However, these atoms still form two main types: the first type consists only of atoms $A$, whilst the second type consists of all other atoms (which are in fact saddle atoms).

Let now $(M, \omega, (H,f))$ be an integrable system on a symplectic 4-manifold $M$ and let $Q = \{x \in M \mid H(x) = \text{constant}\}$ be a `generic' so-called isoenergy $3$-surface (see Bolsinov and Fomenko \cite[Section 3.8]{Bolsinov2004} for the exact conditions on $Q$). Let $Q/\sim$ be the space of leaves, which can also be pictured as a (Reeb) graph where the vertices correspond to the singular leaves. Now, the singular leaves correspond to so-called \emph{$3$-atoms}, which are defined similarly to the atoms we saw before, but now the neighbourhoods are $3$-dimensional. It turns out that these 3-atoms are in one-to-one correspondence with the set of 2-atoms possibly endowed with a finite number of marked points or \textit{stars} -- corresponding to exceptional fibres of the Seifert fibration naturally associated to a 3-atom, see Bolsinov and Fomenko \cite[Sections 3.2 and 3.5]{Bolsinov2004}. For simplicity, 2-atoms with stars will also be referred to as $2$-atoms. Thus, we will consider the graph defined by $Q/\sim$ with the vertices corresponding to $2$-atoms.
This graph is called the \emph{molecule} of $(M, \omega, (H,f))$ on $Q$.

A molecule contains a lot of information of the foliation of the isoenergy surface $Q$. But this type of molecule consists of atoms glued together {\em so far without} the knowledge of how this gluing is performed. Keeping track of the gluing gives us the final piece of information that we need to give a molecule the meaning of an invariant: the gluing is performed by the so-called \emph{gluing matrix}
\begin{align*}
    C_{i} =
    \begin{pmatrix}
        \alpha_{i} & \beta_{i} \\
        \gamma_{i} & \delta_{i}
    \end{pmatrix} \in \textup{GL}(2, \mathbb Z), \quad \det C = -1.
\end{align*}
To the gluing matrix $C_{i}$, there are two invariants assigned, namely
\begin{align*}
    r_{i} := 
    \begin{cases}
        \frac{\alpha_{i}}{\beta_{i}} \mod 1 & \mbox{ if } \beta_{i} \neq 0, \\
        \infty & \mbox{ if }\beta_{i} = 0
    \end{cases}
    \qquad \text{and} \qquad
    \epsilon_{i} := 
    \begin{cases}
        \text{sign}\, \beta_{i} & \mbox{ if } \beta_{i} \neq 0, \\
        \text{sign}\, \alpha_{i} & \mbox{ if } \beta_{i} = 0.
    \end{cases}
\end{align*}
These two invariants alone are not enough for our purposes, and so one more invariant has to be introduced. An edge $e_{i}$ of a molecule $W$ is called \emph{infinite}, if $r_{i} = \infty$, and otherwise \emph{finite}. Cutting the molecule along finite edges splits it into several connected components. The components not containing any atoms of type $A$ are called \emph{families}. Let $U_{k}$ be a family. Recall that the edges of atoms are `oriented' by arrows. An edge in $U_{k}$ is said to be \emph{outgoing} if the arrow points from a vertex inside $U_{k}$ to a vertex outside $U_{k}$. In the opposite case an edge in $U_{k}$ is called \emph{incoming}. If the edge joins a vertex inside $U_{k}$ to another vertex inside $U_{k}$, then the edge is called \emph{interior}. To each edge $e_{i}$ in $U_{k}$ we assign the following integer:
\begin{align*}
    \Theta_{i}: =
    \begin{cases}
        \lfloor \alpha_{i}/\beta_{i} \rfloor, & \mbox{if } e_{i} \text{ is an outgoing edge}, \\
        \lfloor-\delta_{i}/\beta_{i} \rfloor, & \mbox{if } e_{i} \text{ is an incoming edge}, \\
        -\gamma_{i}/\alpha_{i}, & \mbox{if } e_{i} \text{ is an interior edge}.
    \end{cases}
\end{align*}
With this, we construct the third, and final, invariant we want to associate to $W$, namely
\begin{equation*}
    n_{k} := \sum_{e_{i} \in U_{k}} \Theta_{i} \in \Z.
\end{equation*}
The invariants $r_{i}$, $\epsilon_{i}$ and $n_{k}$ will be called \emph{marks}. One can now endow the molecule $W$ with the three marks defined above, and define the \emph{marked molecule} as the quadruple $W^{*} := (W, r_{i}, \epsilon_{i}, n_{k})$. Fomenko and Zieschang \cite{Fomenko1990} showed that two integrable systems on generic isoenergy 3-surfaces are Liouville equivalent if and only if their marked molecules coincide. Marked molecules are also known as \emph{Fomenko-Zieschang invariants}. The collection of such marked molecules can be thought of as a topological portrait of the system, which contains more information than for example the topological types of the individual singular leaves/fibres. 

Since hypersemitoric systems only contain elliptic, hyperbolic-regular, focus-focus and parabolic points, but no hyperbolic-hyperbolic ones, one can show that marked loop molecules form complete local topological invariants of the torus fibration of a hypersemitoric system. In other words, the loop molecules around a given 
singularity of the hypersemitotic system determine its topological type.
Note that the same is not true for general hyperbolic-hyperbolic singularities of
integrable 2 degree of freedom systems; see Bolsinov and Oshemkov \cite{Bolsinov2006}.
\section{Semitoric systems}\label{sec:semitoric}
\noindent
In this section, we will briefly recall the construction of the five invariants of semitoric systems introduced by Pelayo and V{\~u} Ng{\d{o}}c \cite{Pelayo2009} and its generalizations, then observe transitions from toric to semitoric systems by creating focus-focus points, and eventually consider some explicit examples.

Two semitoric systems $(M_{1},\omega_{1},(J_{1},H_{1}))$ and $(M_{2},\omega_{2},(J_{2},H_{2}))$ are said to be \emph{isomorphic} if there exists a symplectomorphism $\varphi : M_{1} \to M_{2}$ such that $\varphi^{*}(J_{2},H_{2}) = (J_{1},f(J_{1},H_{1}))$ for some smooth function $f$ such that $\frac{\partial f}{\partial H_{1}} > 0$. Since semitoric systems always come with a smooth, globally defined action $J$, this definition is basically saying that two semitoric systems are equivalent if and only if the corresponding Lagrangian fibrations are fibrewise symplectomorphic (up to possibly changing $J$ to $\pm J + \textup{const}$).  

Pelayo and V{\~u} Ng{\d{o}}c \cite[Theorem 6.2]{Pelayo2009} showed that two simple semitoric systems are isomorphic if and only if all five invariants (defined below) are equal for the two systems. The simplicity assumption has been removed from the classification by Palmer, Pelayo and Tang \cite{Palmer2019}, but the invariants in the non-simple case are more complicated, and we do not present them here.

\subsection{The five semitoric invariants}

Let $(M, \omega, F=(J,H))$ be a simple semitoric system. We will use the identification $S^{1} = \R/2\pi\Z$ in what follows.
Let us now explain each of the five invariants in more detail.

\subsubsection{Number of focus-focus points}

V{\~u} Ng{\d{o}}c \cite[Corollary 5.10]{Vu-Ngoc2007} proved that $M$ has a finite number of focus-focus singular points. Denoting this number by $n_{FF}$, one has thus $0 \leq n_{FF} < \infty$. Then $n_{FF}$ forms an invariant for semitoric systems (cf.\ Pelayo and V{\~u} Ng{\d{o}}c \cite[Lemma 3.2]{Pelayo2009}). 

\subsubsection{Taylor series invariant} \label{sec:taylor_series_inv}

Denote the focus-focus points of $(M, \omega, F=(J,H))$ by $m_i$ for $1 \leq i \leq n_{FF}$. Let us now consider one focus-focus point, and denote it by $m$ without the index, to simplify the notation. Recall from Section \ref{sec:IntSys} that there exists a neighbourhood $U$ of $m$ with symplectic coordinates $(x,y,\xi,\eta)$ such that the quadratic parts of $J$ and $H$ span a Cartan subalgebra with the following basis:
\begin{equation} \label{eq:local-focus-focus}
    q_{1} = x\xi + y\eta, \qquad q_{2} = x\eta - y\xi.
\end{equation}
Note that the Hamiltonian flow generated by $q_{2}$ is $2\pi$-periodic.

We now follow the exposition in V{\~u} Ng{\d{o}}c \cite{Vu-Ngoc2003}:
Let $\Lambda_{z} = F^{-1}(z)$ be a regular fibre near the singular fibre containing $m$. For any point $A \in \Lambda_{z}$, denote by $\tau_{1}(z)$ the first return time of the flow generated by $X_{H}$ to the $X_{J}$-orbit through $A$, and let $\tau_{2}(z) \in \R/2\pi\Z$ be the time it takes to close up this trajectory under the flow of $X_{J}$. V{\~u} Ng{\d{o}}c \cite[Proposition 3.1]{Vu-Ngoc2003} showed that, for some determination of the complex logarithm $\ln z$, then 
\begin{equation} \label{eq:sigma}
    \sigma_{1}(z) := \tau_{1}(z) + \Re(\ln z), \quad
    \sigma_{2}(z) := \tau_{2}(z) - \Im(\ln z)
\end{equation}
extends to smooth and single-valued functions in a neighbourhood of $c = F(m)$. Moreover, $\sigma := \sigma_{1} dz_{1} + \sigma_{2} dz_{2}$ yields a closed $1$-form under the identification $z=(z_1, z_2) \in \R^2$. 
Define $S$ via $dS = \sigma$ and $S(c) = 0$ and denote the Taylor series of $S$ at $z = c$ by $(S)^{\infty}$. The {\em Taylor series invariant}, for all focus-focus points $m_{i}$, $1 \leq i \leq n_{FF}$, is then given by the $n_{FF}$-tuple $((S_{i})^{\infty})_{i=1}^{n_{FF}}$.

%Note that $\tau_{i}^{1}(z)$ shows a logarithmic divergence when $z$ approaches the focus-focus value $c_i$. 
 There is another way to define the Taylor series invariant. Let $\gamma_{z}^1$ and $\gamma_{z}^2$ be a basis of the first homology group of the torus $\Lambda_z$ that varies smoothly with the base point $z$ such that $\gamma_{z}^1$ is a representative of the cycle corresponding to the (periodic) flow of $J$ and $\gamma_{z}^2$ represents a homology cycle obtained by first moving with the flow of $X_H$ using time  
$\tau_{1}(z)$ and then with the flow of $X_J$ using time $\tau_{2}(z)$.
Now consider the action integral
\begin{align*}
    \mathcal{A}(z) := \int_{\gamma_{z}^2} \alpha,
\end{align*}
where $\alpha$ is a primitive of $\omega$ on some neighbourhood of $\Lambda_{z}$. Then one finds for $z\simeq(z_1,z_2) \in \R^2$
\begin{equation} \label{eq:dAction}
    d\mathcal{A}(z) = \tau_{1}(z) dz_{1} + \tau_{2}(z) dz_{2}.
\end{equation}
One can in fact interpret $S$ as a \emph{regularised action integral} via
\begin{equation} \label{eq:regularised_action}
    S(z) = \mathcal{A}(z) - \mathcal{A}(c) + \Re(z \ln z - z).
\end{equation}
Note that the above construction involves a certain number of choices which have to be made compatibly with the construction of the polygon invariant and the twisting index invariant below. The exact dependencies are explained in detail in the forthcoming article by Alonso, Hohloch, and Palmer \cite{AlonsoHohlochPalmer}.

\subsubsection{Polygon invariant}
\label{sec:polygonInvariant}
Let $m_1, \dots, m_{n_{FF}}$ be the focus-focus points and denote by $c_1:=F(m_1)$, \dots, $c_{n_{FF}}:= F(m_{n_{FF}})$ their values ordered such that the first coordinate of the focus-focus values increases. Denote by $B := F(M)$ the image of the momentum map. V{\~u} Ng{\d{o}}c \cite[Theorem 3.4]{Vu-Ngoc2007} showed that the set $B_r \subseteq F(M)$ of regular values of $F$ coincides with the set $ \textnormal{int}\, B \setminus \{c_{1}, \dots, c_{n_{FF}}\}$. One can render $B_{r}$ simply connected by making a vertical cut from each focus-focus value $c_{i}$ either upwards or downwards to the boundary of $F(M)$.

By the Arnol'd-Liouville theorem, the momentum map induces an integral affine structure on $B$ (which in general does not agree with the one induced by the inclusion of $B$ into $\R^2$). Recall that affine transformations leaving a vertical line invariant arise from vertical translations composed with a matrix of the form
\begin{align*}
    T^{k} := \begin{pmatrix}
        1 & 0 \\
        k & 1
    \end{pmatrix}
\end{align*}
with $k \in \Z$. Now denote by $l_{i} \subset \R^{2}$ the vertical line through the focus-focus singular value $c_{i} \in \R^2$. This line splits $\R^{2}$ into two half-spaces. For $k \in \Z$, let $t_{l_{i}}^{k} : \R^{2} \to \R^{2}$ be the map that leaves the left half-space invariant and shears the right half-space by $T^{k}$. We accommodate now all focus-focus singular values by setting $\mathbf{k} := (k_{1}, \dots, k_{n_{FF}})$ and defining $t_{\mathbf{k}} := t_{l_{1}}^{k_{1}} \circ \dots \circ t_{l_{n_{FF}}}^{k_{n_{FF}}}$. 

For each $1 \leq i \leq n_{FF}$, let $\epsilon_{i} \in \{-1, +1\}$, and denote by $l_{i}^{\epsilon_{i}}$ the vertical half line starting at $c_{i}$, going upwards if $\epsilon_{i} = +1$, and downwards if $\epsilon_{i} = -1$, and let $l^{\epsilon} := l_1^{\epsilon_1} \cup\ \dots\ \cup\ l_{n_{FF}}^{\epsilon_{n_{FF}}}$ be the union of the lines running through all focus-focus values for a choice of $\epsilon := (\epsilon_{1}, \dots , \epsilon_{n_{FF}})$. Then the set $B \setminus l^\epsilon$ is simply connected for all possible choices of $\epsilon_i$. 

V{\~u} Ng{\d{o}}c \cite[Theorem 3.8]{Vu-Ngoc2007} showed that there exists a homeomorphism $f:=f_\epsilon : B \to \R^{2}$ depending on the choices of $\epsilon$ and preserving $J$ such that $f(B)$ is a rational convex polygon. 
Restricted to $B\setminus l^\epsilon$, the homeomorphism $f$ becomes a diffeomorphism onto its image which sends the integral affine structure of $B_r \setminus l^\epsilon$ to the integral affine structure of $\R^2$.
The map $\mu := f \circ F$ is called a \emph{generalized toric momentum map} for $(M, \omega, F=(J,H))$ (cf.\ Pelayo and V{\~u} Ng{\d{o}}c \cite[Definition 4.3]{Pelayo2009}). 

In order to turn the polygon $\Delta := \mu(M)$ into an invariant of the underlying semitoric system one needs to get rid of the choices involved in the construction of $\Delta$. This is done by means of a group action: consider the group $\mathcal{G} := \{T^{k} \mid k \in \Z\}$ and the action of the group $\{-1, +1\}^{n_{FF}} \times \mathcal{G}$ on $\bigl(\Delta, (l_{i})_{i=1}^{n_{FF}}, (\epsilon_{i})_{i=1}^{n_{FF}}\bigr)$ given by
\begin{align*}
    \bigl((\epsilon'_{i})_{i=1}^{n_{FF}}, T^{k} \bigr) \cdot \bigl(\Delta, (l_{i})_{i=1}^{n_{FF}}, (\epsilon_{i})_{i=1}^{n_{FF}} \bigr) \
    := \ \bigl(t_{\mathbf{u}}(T^{k}(\Delta)), (l_{i})_{i=1}^{n_{FF}}, (\epsilon'_{i}\epsilon_{i})_{i=1}^{n_{FF}}\bigr)
\end{align*}
where $\mathbf{u} = ((\epsilon_{i}- \epsilon'_{i})/2)_{i=1}^{n_{FF}}$. Then the {\em polygon invariant} is the orbit of $\bigl(\Delta, (l_{i})_{i=1}^{n_{FF}}, (\epsilon_{i})_{i=1}^{n_{FF}}\bigr)$ under the above action (cf.\ Pelayo and V{\~u} Ng{\d{o}}c \cite[Definition 4.5 and Lemma 4.6]{Pelayo2009}).

\subsubsection{Height invariant}

For $i \in \{1, \dots, n_{FF}\}$, consider the focus-focus singular points $m_i$ and their images $c_{i} := F(m_{i})$ and let $\mu$ and $\Delta$ be as in Section \ref{sec:polygonInvariant}. The \textit{height} (or the \textit{volume}) invariant, as introduced by Pelayo and V{\~u} Ng{\d{o}}c \cite[Definition 5.2 and Lemmas 5.1 and 5.3]{Pelayo2009}, is given by the $n_{FF}$-tuple $(h_1, \dots, h_{n_{FF}})$ with
\begin{equation*}
    h_{i} := \text{pr}_2(\mu(m_{i})) - \min_{s \in l_{i} \cap \Delta} \text{pr}_{2}(s),
\end{equation*}
where $\text{pr}_{2} : \R^{2} \to \R$ is the projection onto the second coordinate (in \cite[Remark 5.2]{Pelayo2009} it is explained how this height invariant corresponds to the volume of certain submanifolds, and hence it is sometimes called the volume invariant). The function $h_i$ thus measures the distance between the focus-focus value in the polygon $\Delta=\mu(M)$ and its lower boundary. Furthermore, $h_i$ is independent of the choice of the generalized toric momentum map $\mu$, since it can also be seen as the symplectic volume of certain level sets.

 \subsubsection{Twisting index invariant}
\iffalse
Let $U_{i}$ be a neighbourhood of a focus-focus singular point $m_{i} \in F^{-1}(c_{i})$, and let $V_{i} = F(U_{i})$. V{\~u} Ng{\d{o}}c and Wacheux \cite{Vu-Ngoc2013} showed that there exists a local symplectomorphism $\Psi : (\R^{4}, \omega_{0}) \to (M, \omega)$ sending the origin to $m_{i}$, and a local diffeomorphism $G : \R^{2} \to \R^{2}$ sending $0$ to $F(m_{i})$ such that $F \circ \Psi = G \circ q_{i}$, where $q_{i} = (q_{i}^{1}, q_{i}^{2})$ is given by \eqref{eq:local-focus-focus}. Recall that $q_{i}^{2}$ generates a circle action, so it must correspond to $J$. If necessary, after composing $\Psi$ with either/both of the canonical transformations $(x, \xi) \mapsto (-x, -\xi)$ and $(x, y, \xi, \eta) \mapsto (-\xi, -\eta, x, y)$, one finds that $G$ is of the form
\begin{equation*}
    G(q_{i}^{1}, q_{i}^{2}) = (q_{i}^{2}, G_{2}(q_{i}^{1}, q_{i}^{2})),
\end{equation*}
where $\del[G_{2}]{q_{i}^{1}}(0) > 0$. We will extend $G_{2}(q_{i}^{1}, q_{i}^{2})$ to another Hamiltonian function $\Tilde{G}_{2}(\Tilde{H}, J)$, such that they are equal at their restriction to $U_{i}$. Here $(\Tilde{H}, J)$ is a new momentum map for the semitoric system, and $\Tilde{G}_{2} : \R^{2} \to \R$ is some function to be discussed further below.
\fi

Recall the action integral introduced in the construction of the Taylor series invariant (see Subsection~\ref{sec:taylor_series_inv}):
\begin{equation*}
   \mathcal{A}_{i}(z) := \int_{\gamma_{i, z}^2} \alpha.
\end{equation*}
Let $G_{i}(z) := \mathcal{A}_{i}(z) - \mathcal{A}_{i}(c_{i})$ for $ i = 1, \ldots, n_{FF}$. 
Observe that $G_{i}(0)$ is well defined and equal to zero since the actions $\mathcal{A}_{i}(z)$ are given by integrating a primitive 1-form over a loop on a Lagrangian torus $\Lambda_{z}$. Note that this could also have been seen by using the regularised action in \eqref{eq:regularised_action}. Now, let us define the Hamiltonian function via $H_{i, p} := G_i(J, H).$ Then $\lim_{m \to m_{i}} H_{i, p} = 0$. Note also that, by \eqref{eq:dAction}, we get a Hamiltonian vector field
\begin{equation*}
    X_{i, p} = (\tau_{i}^{1} \circ F) X_{J} + (\tau_{i}^{2} \circ F) X_{H}.
\end{equation*}
This was discussed by Pelayo and V{\~u} Ng{\d{o}}c \cite[Section 5.2]{Pelayo2009}. They called the momentum map $\nu := (J, H_{i, p})$ the \emph{privileged momentum map} for $F = (J, H)$.

Now, let $\mu$ be a generalized toric momentum map. As $\mu$ preserves $J$, its components satisfy $(\mu_{1}, \mu_{2}) = (J, \mu_{2})$. As $\mu_i, J$ and $H_{i,p}$ are all action variables, there exists an invertible matrix $A \in \textnormal{GL}(2, \Z)$ such that $(X_{J}, X_{\mu_{2}}) = A(X_{J}, X_{i, p})$. The matrix has to be of the form
\begin{equation*}
    A = 
    \begin{pmatrix}
    1 & 0 \\
    k_{i} & 1
    \end{pmatrix},
\end{equation*}
hence $X_{\mu_{2}} = k_{i} X_{J} + X_{i, p}.$
Pelayo and V{\~u} Ng{\d{o}}c \cite[Proposition 5.4]{Pelayo2009} showed that $k_{i}$ does not depend on $X_{i, p}$ or $G_i$. The integer $k_{i}$ is called the \emph{twisting index}. Note that, if $k_{i}$ is the twisting index of $m_{i}$, then locally $\mu = T^{k_{i}} \nu$. Also, if the polygon is transformed by some $T^{r}$, then $\nu$ does not change, whilst $\mu \to T^{r} \mu$.

Note that the twisting index depends on the polygon $\Delta$. To introduce an actual invariant, similarly to Subsection~\ref{sec:polygonInvariant}, 
we consider the orbit of $(\Delta, (l_{i})_{i=1}^{n_{FF}}, (\epsilon_{i})_{i=1}^{n_{FF}}, (k_{i})_{i=1}^{n_{FF}})$ under the action of $\{-1, +1\}^{n_{FF}} \times \mathcal{G}$. Specifically, with $\mathbf{u} := (u_i)_{i=1}^{n_{FF}} := ((\epsilon_{i}-\epsilon_{i}\epsilon'_{i})/2)_{i=1}^{n_{FF}}$, the action is given by
\begin{align*}
    & ((\epsilon'_{i})_{i=1}^{n_{FF}}, T^{k}) \cdot (\Delta, (l_{i})_{i=1}^{n_{FF}}, (\epsilon_{i})_{i=1}^{n_{FF}}, (k_{i})_{i=1}^{n_{FF}}) \\
    & \quad \qquad \qquad \qquad \qquad = \left(t_{\mathbf{u}}(T^{k}(\Delta)), (l_{i})_{i=1}^{n_{FF}}, (\epsilon'_{i} \epsilon_{i})_{i=1}^{n_{FF}}, \left(k + k_{i} + \sum_{j=1}^{\tilde{i}_i} u_j\right)_{i=1}^{n_{FF}} \right)
\end{align*}
where we set $0=:\sum_{j=1}^{0} u_j$ and where $\tilde{i}_i =i $ or $\tilde{i}_i= i-1$ depending on the choice of certain conventions.
This orbit is called the {\em twisting index invariant} (cf.\ Pelayo and V{\~u} Ng{\d{o}}c \cite[Definition 5.9 and Lemma 5.10]{Pelayo2009}).
Note that the above formula differs slightly from the original one given in Pelayo and V{\~u} Ng{\d{o}}c \cite{Pelayo2009} by the extra term $\sum_{j=1}^{\tilde{i}_i} u_j$. This term accounts for the way in which changing cut directions affects the twisting index.
Its absence in the original formula was pointed out to us by Yohann Le Floch and Joseph Palmer  (for a detailed discussion, we refer to the forthcoming paper by Alonso, Hohloch, and Palmer \cite{AlonsoHohlochPalmer}).

\subsection{Modifications and generalizations of the five invariants}

In fact, all five invariants are intimately related, and there is no need to consider them separately. Le Floch and Palmer \cite{LeFloch2018} took three of the five semitoric invariants --- the number of focus-focus points, the polygon invariant, and the height invariant --- and joined them together to form a single invariant, called the \emph{marked semitoric polygon invariant}. When Palmer, Pelayo and Tang \cite{Palmer2019} extended the classification to non-simple semitoric systems they gathered all five invariants into one big invariant, called the \emph{complete semitoric invariant}.

\subsection{Supercritical Hamiltonian-Hopf bifurcation}\label{sec:supercritical_Hamiltonian-Hopf}

If one perturbs a toric system, one may obtain a semitoric system, in particular if an elliptic-elliptic point is transformed into a focus-focus point. Such a transformation is called a \emph{supercritical Hamiltonian-Hopf bifurcation}. In 
coordinate form, it can more specifically be defined as follows
(see in particular Equation \eqref{eq:Hamiltonian-Hopf} below with $a >0$).

Let $\mathfrak G$ be a Lie group acting on the space of smooth real-valued functions $C^{\infty}(\R^{n})$ whose action is defined by $g \cdot f(x) = f(g^{-1}(x))$ for $g \in \mathfrak  G$, $f \in C^{\infty}(\R^{n})$ and $x \in \R^{n}$. Furthermore, let $\R[x]$ denote the space of polynomials on $\R^{n}$, and let $\R[x]^{\mathfrak G}$ be the space of $\mathfrak G$-invariant polynomials. Hilbert showed that, if ${\mathfrak G}$ is compact, then there exist finitely many invariant polynomials $\rho_{i} \in \R[x]^{\mathfrak  G}$ for $i = 1, \dots, k$ which generate $\R[x]^{\mathfrak G}$ as an algebra (cf.\ van der Meer \cite[Section 3.1]{vanderMeer1985}). Such invariant polynomials $\rho_{i}$ are called \emph{Hilbert generators}.

Let $(x, y, \xi, \eta)$ be canonical coordinates on $\R^{4}$ and define the following three Hilbert generators: $J = x \eta - y \xi$, $X = \frac{1}{2}(\xi^{2} + \eta^{2})$, and $Y = \frac{1}{2}(x^{2} + y^{2})$. When considering (hyper)semitoric systems, we will choose $\mathfrak G = S^{1}$ to be given by the periodic Hamiltonian flow of $X_J$.  Then van der Meer \cite[Corollary 3.39]{vanderMeer1985} showed that there exists the following equivariant normal form for a Hamiltonian-Hopf bifurcation
\begin{equation}\label{eq:Hamiltonian-Hopf}
    \hat{H}_{s} = J + X + s Y + a Y^{2},
\end{equation}
where $s, a \in \R$ are parameters with $a \neq 0$, which we for simplicity take as a definition for this type of bifurcation. 
If $a > 0$ the bifurcation is called \textit{supercritical}, and
\textit{subcritical} otherwise. Note that here the momentum map is given by $(J, \hat H_s).$

Recall that the singular points in a $2$-degree of freedom toric system all have only elliptic and/or regular components. If we perturb one of the integrals of a $2$-degree of freedom toric system as in the above normal form, then we can make one of the elliptic-elliptic singular points turn into a focus-focus point. On the level of eigenvalues, 4 purely imaginary eigenvalues at an elliptic-elliptic point collide when the bifurcation parameter attains the value $s = 0$ and then
change into four complex eigenvalues (cf.\ van der Meer \cite[Section 1.3]{vanderMeer1985}). One can see two examples of supercritical Hamiltonian-Hopf bifurcations in Figure \ref{fig:cam_moment_image} and Figure \ref{fig:octagon_polygon}. The subcritical case, when the sign of $a$ is negative, is treated in Section \ref{sec:subcritical_Hamiltonian-Hopf}.

\subsection{Examples}

To compute the semitoric invariants explicitly for given systems has proven to be very difficult since it needs the combination of theoretical knowledge and strong computational skills.

\subsubsection{Coupled angular momenta system}\label{sec:semitoric_cam}

Consider the manifold $M := S^{2} \times S^{2}$ and equip it with the symplectic form $\omega := - (R_{1}\omega_{S^{2}} \oplus R_{2}\omega_{S^{2}})$ where $\omega_{S^{2}}$ is the standard symplectic form on $S^{2}$ and $R_{1}, R_{2} \in \R^{>0}$. 
When Sadovski{\'\i} and Zhilinski{\'\i} \cite{Sadovskii1999} studied the so-called {\em coupled angular momenta system}, they found a focus-focus point and nontrivial monodromy. Since this system is both interesting from a physics point of view and not very complicated from a mathematical point of view, it recently became a popular subject to study.

\begin{figure}[h]
    \subfloat{\includegraphics[scale=0.3]{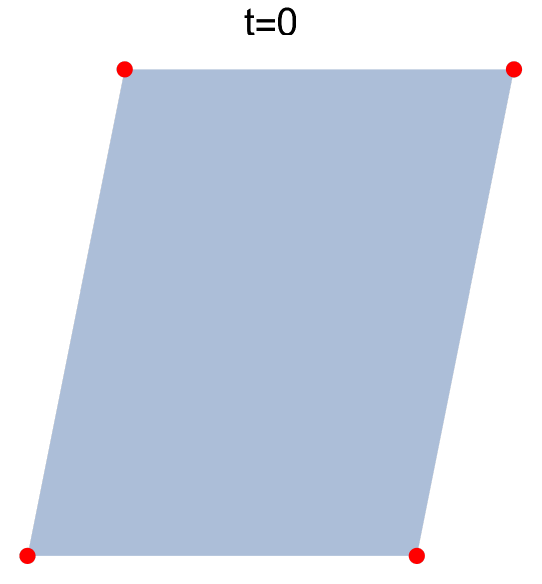}}
    \subfloat{\includegraphics[scale=0.3]{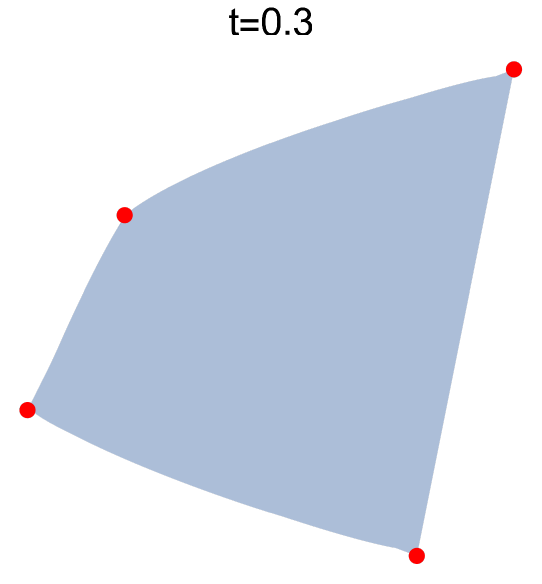}}
    \subfloat{\includegraphics[scale=0.3]{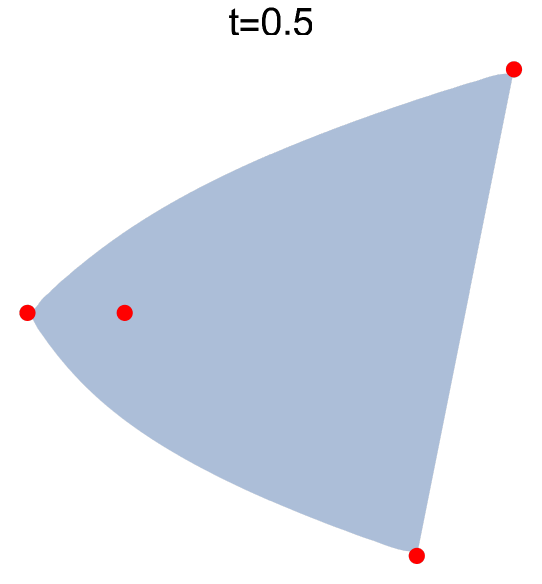}}
    \subfloat{\includegraphics[scale=0.3]{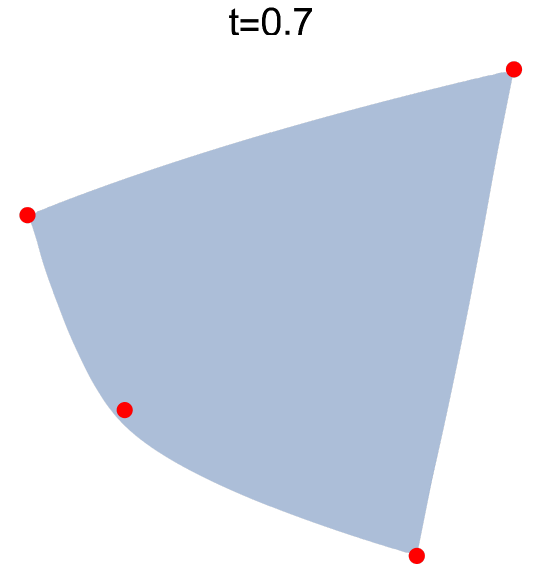}}
    \subfloat{\includegraphics[scale=0.3]{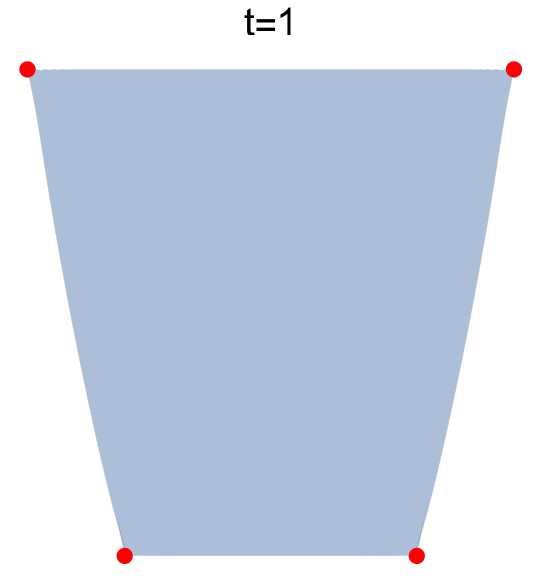}}
    \caption{The momentum map image for the coupled angular momenta with $R_{1} = 1$ and $R_{2} = 4$, for selected values of the coupling parameter $t$. The red points mark the images of the four rank zero points, which are of elliptic-elliptic type unless they are located in the interior of the polygon (the blue area) where they are of focus-focus type. When passing from the second to the third subfigure from the left, a supercritical Hamiltonian-Hopf bifurcation takes place.}
    \label{fig:cam_moment_image}
\end{figure}

Le Floch and Pelayo \cite{LeFlochPelayoCAM} showed that the coupled angular momenta system on $M$, given in Cartesian coordinates by
\begin{align*}
    J(x_{1},y_{1},z_{1},x_{2},y_{2},z_{2}) &:= R_{1}(z_{1}-1) + R_{2}(z_{2}+1), \\
    H(x_{1},y_{1},z_{1},x_{2},y_{2},z_{2}) &:= (1-t)z_{1} + t(x_{1}x_{2} + y_{1}y_{2} + z_{1}z_{2}),
\end{align*}
describes a semitoric system for all $t \in \R \setminus \{t^{-},t^{+}\}$, where
\begin{equation*}
    t^{\pm} := \frac{R_{2}}{2R_{2} + R_{1} \mp 2\sqrt{R_{1}R_{2}}}.
\end{equation*}
The system has four singular points of rank 0 which are located at the top and bottom of the spheres, i.e.\ when $(z_{1},z_{2}) = (\pm 1, \pm 1)$. Three of the points are always elliptic-elliptic, whilst $(1, -1)$ is a focus-focus point if $t^{-} < t < t^{+}$ and elliptic-elliptic if $t < t^{-}$ or $t > t^{+}$. Thus, the number of focus-focus points invariant is $0$ if $(1, -1)$ is elliptic-elliptic, or $1$ if $(1, -1)$ is focus-focus.
For some values of $t$, the moment image is plotted in Figure \ref{fig:cam_moment_image}. 

\begin{figure}[h]
    \centering
    \includegraphics[scale=0.8]{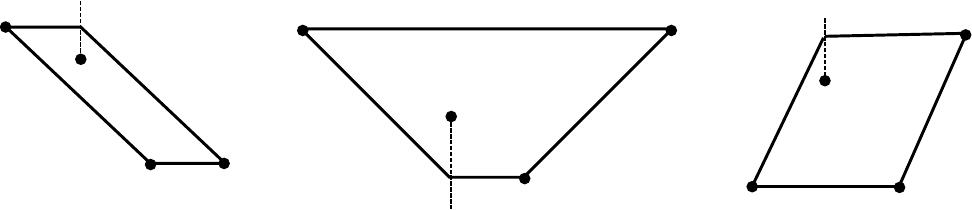}
    \caption{This figure displays some of the polygons of the polygon- and twisting index invariant, for the case $R_{1} < R_{2}$. From left to right, the values for the twisting index is $-1$, $-1$ and $0$. The dashed lines indicates the cutting direction. For plots of all three cases $R_{1} > R_{2}$, $R_{1} = R_{2}$ and $R_{1} = R_{2}$, see Alonso, Dullin and Hohloch \cite[Figure 1]{Alonso2020a}.}
    \label{fig:cam_polygon}
\end{figure}

Le Floch and Pelayo \cite{LeFlochPelayoCAM} computed, for certain parameter values, the first two terms of the Taylor series, the polygon, and the height invariant for this system. The full classification was achieved by Alonso, Dullin and Hohloch \cite{Alonso2020a}. The semitoric invariants of the coupled angular momenta system are as follows: The number of focus-focus points is either zero or one, see above. The Taylor series invariant is of the form
\begin{align*}
    S(j,k) =& 
    j \arctan \left( \frac{R_{2}^{2}(2t - 1) - R_{1}R_{2}(t + 1) + R_{1}^{2}t}{(R_{1} - R_{2})R_{1} r_{A}} \right) 
    + k \ln \left( \frac{4 R_{1}^{5/2} r_{A}^{3}}{R_{2}^{3/2}(1 - t) t^{2}} \right) \\&
    + \frac{j^{2}}{16 R_{1}^{4} R_{2} r_{A}^{3}} \Big( R_{2}^{4}(2t - 1)^{3} - R_{1}R_{2}^{3}(32t^{3} - 46t^{2} + 17t - 1) \\& \hspace{2.6cm}
    - 3R_{1}^{2}R_{2}^{2}t(4t^{2} - 7t + 1) + R_{1}^{3}R_{2}(3 - 5t)^{2} - R_{1}^{4}t_{3} \Big) \\&
    + \frac{jk(R_{2} - R_{1})}{8R_{1}^{3}R_{2}r_{A}^{3}} \left( R_{2}^{2}(2t - 1)^{2} - 2R_{1}R_{2}t(6t - 1) + R_{1}^{2}t^{2} \right) \\&
    + \frac{k^{2}}{16R_{1}^{4}R_{2}r_{A}^{3}} \Big( R_{2}^{4}(2t - 1)^{3} - R_{1}R_{2}^{3}(16t^{3} - 42t^{2} + 15t + 1) \\& \hspace{2.6cm}
    - R_{1}^{2}R_{2}^{2}t(28t^{2} - 3t -3) + R_{1}^{3}R_{2}t^{2}(13t - 3) + R_{1}^{4}t^{3} \Big) \\&
    + \mathcal{O}(3),
\end{align*}
where
\begin{align*}
    r_{A} = \sqrt{(R_{1}^{2} + 4R_{2}^{2})(t - t^{-})(t^{+} - t)}.
\end{align*}
The polygon and twisting index invariants are illustrated in Figure \ref{fig:cam_polygon}. 

Set $R:= \frac{R_2}{R_1}$. Then the height invariant of the coupled angular momenta is given by
\begin{align*}
    h =  2 \min(R_{1}, R_{2}) 
    + \frac{R_{1}}{\pi t} \left( r_{A} - 2 R t \arctan \left( \frac{r_{A}}{R - t} \right) - 2 t \arctan \left( \frac{r_{A}}{R + t - 2 R t} \right) \right).
\end{align*}

\subsubsection{The (semi)toric octagon system} \label{sec:semitoric_octagon}

\begin{figure}[h]
    \centering
    \includegraphics{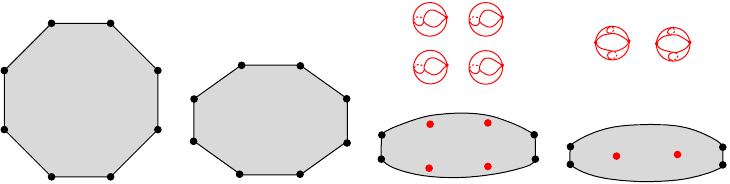}
    \caption{This figure shows the octagon system as it transitions from having eight elliptic-elliptic values in the first two figures, to having four elliptic-elliptic values and four focus-focus values in the third figure, and finally four elliptic-elliptic values and two focus-focus values whose fibres both contains two focus-focus points. Above the final two figures, the fibres of the focus-focus values are drawn.}
    \label{fig:octagon_polygon}
\end{figure}

De Meulenaere and Hohloch \cite{Meulenaere2019} constructed a semitoric system with four focus-focus singular points. The system was created by first considering the octagon $\Delta$ obtained by chopping off the corners of the square $[0, 3] \times [0,3]$. Since $\Delta$ turned out to be a Delzant polygon, Delzant's \cite{Delzant1988} construction could be used to construct a toric system which has $\Delta$ as image of the momentum map. This is done by means of symplectic reduction of $\C^8$ (equipped with its standard symplectic structure) and yields a $4$-dimensional, compact, connected, symplectic manifold $(M_{\Delta}, \omega_{\Delta})$. A point on $M_{\Delta}$ is written as an equivalence class of the form $[z] = [z_{1}, \dots, z_{8}]$ with $z_{i} \in \C$ for $i = 1, \dots, 8$. The (toric) momentum map $F = (J, H):(M_{\Delta}, \omega_{\Delta}) \to \R^2 $ is given by
\begin{align*}
    J([z_{1}, \dots, z_{8}]) = \frac{1}{2}\abs{z_{1}}^{2}, \quad
    H([z_{1}, \dots, z_{8}]) = \frac{1}{2}\abs{z_{3}}^{2}.
\end{align*}
Denote by $\Re $ the real part of a complex number. By perturbing $H$ to 
$$H_{t}: = (1-2t) H + t \gamma \Re\left( \bar{z}_{2}\bar{z}_{3}\bar{z}_{4}z_{6}z_{7}z_{8} \right)$$
for $0 < \gamma < \frac{1}{48}$, De Meulenaere and Hohloch \cite[Theorem 4.7]{Meulenaere2019} obtained a system with momentum map $(J, H_{t}):(M_{\Delta}, \omega_{\Delta}) \to \R^2$ that is toric for $0 \leq t < t^{-}$, semitoric for $t^{-} < t < t^{+}$, and toric again for $t^{+} < t \leq 1$, where
\begin{align*}
    t^{-} := \frac{1}{2(1 + 24 \gamma)} \quad \text{and} \quad 
    t^{+} := \frac{1}{2(1 - 24 \gamma)}.
\end{align*}
Note that $0 < t^{-} < \frac{1}{2}$ and $\frac{1}{2} < t^{+} < 1$. At $t = \frac{1}{2}$, the system has two focus-focus fibres, each containing two focus-focus points, see Figure \ref{fig:octagon_polygon}. The two fibres then have the shape of double pinched tori. Apart from one representative of the polygon invariant and the number of focus-focus point, no semitoric invariants have yet been calculated.

\subsection{State of the art concerning other semitoric systems}

Spread over the literature (cf.\ works by Babelon, Dullin, Le Floch, Pelayo, V{\~u} Ng{\d{o}}c, and others), there are various partial results concerning the computation of the semitoric invariants for certain parameter values for certain systems.

For instance, a Taylor series type invariant has been calculated by Dullin \cite{DullinSpherPendulum} for the {\em spherical pendulum} (which is, strictly speaking, not a semitoric system due to lack of properness). 

Pelayo and V{\~u} Ng{\d{o}}c \cite{PelayoVuNgocCSO} computed the number of focus-focus points, the polygon, and the height invariant for the so-called {\em coupled spin oscillator} system. Alonso, Dullin and Hohloch \cite{Alonso2019a} completed the set of semitoric invariants for this system by computing the Taylor series and twisting index invariant.

Both of these systems have only one focus-focus point. Hohloch and Palmer \cite{HohlochPalmer2018} generalized the coupled angular momenta system to a family of semitoric systems with two focus-focus points. Alonso and Hohloch \cite{Alonso2020b} computed the polygon and height invariant for a subfamily and Alonso, Hohloch and Palmer \cite{AlonsoHohlochPalmer} are currently computing its twisting index invariant.

Le Floch and Palmer \cite{LeFloch2018} devised semitoric systems arising from Hirzebruch surfaces and computed their number of focus-focus points, the polygon invariant, and, for certain parameter values, also their height invariant.

\section{Hypersemitoric systems}\label{sec:hypersemitoric}

In this section, we give a brief overview of existing and related results concerning hypersemitoric systems. Recall that, compared to semitoric systems, a hypersemitoric system (Definition \ref{def:HypersemitoricSystem}) may in addition have singular points with hyperbolic components and degenerate singular points of parabolic type.

\subsection{Flaps and pleats/swallowtails} \label{sec:flaps_pleats}

Two possibilities of how hyperbolic-regular and parabolic points occur in hypersemitoric systems are so-called {\em flaps} and {\em pleats/swallowtails}. A good exposition with examples for pleats/swallowtails can be found in Efstathiou and Sugny \cite{EfstathiouSugny}, and for flaps see Efstathiou and Giacobbe \cite{Efstathiou2012}.

There are various ways to visualize flaps and pleats/swallowtails. Instead of using the image of the momentum map over which a hypersemitoric (or even more general) system gives rise to a singular fibration with possibly disconnected fibres, it makes sense to remember the branching and disconnectedness by working with the so-called {\em bifurcation complex} (also known as {\em unfolded momentum domain}). 
One can either identify it with the leaf space of a system $(M, \omega, F=(J, H))$ or describe it directly as a stratified manifold $V$ together with a map $\tilde{F}: M \to V$ and a projection $\tau: V \to \R^2$ such that $\tau \circ \tilde{F} = F$ and the regular level sets of $\tilde{F}$ correspond to the connected components of the level sets of $F$. We will summarize some of their findings.

\begin{figure}[ht]
    \centering
    \begin{subfigure}[t]{0.3\linewidth}
    \includegraphics[scale=0.3]{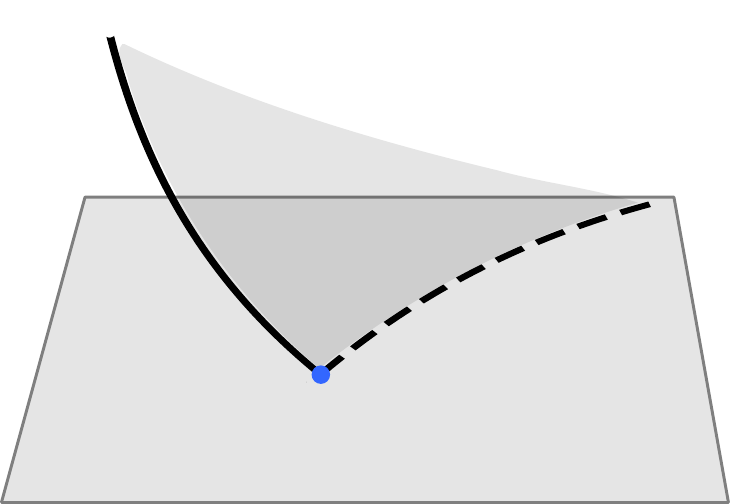}
    \caption{A local flap.}
    \label{fig:local_flap}
    \end{subfigure}
    \begin{subfigure}[t]{0.3\linewidth}
    \includegraphics[scale=0.3]{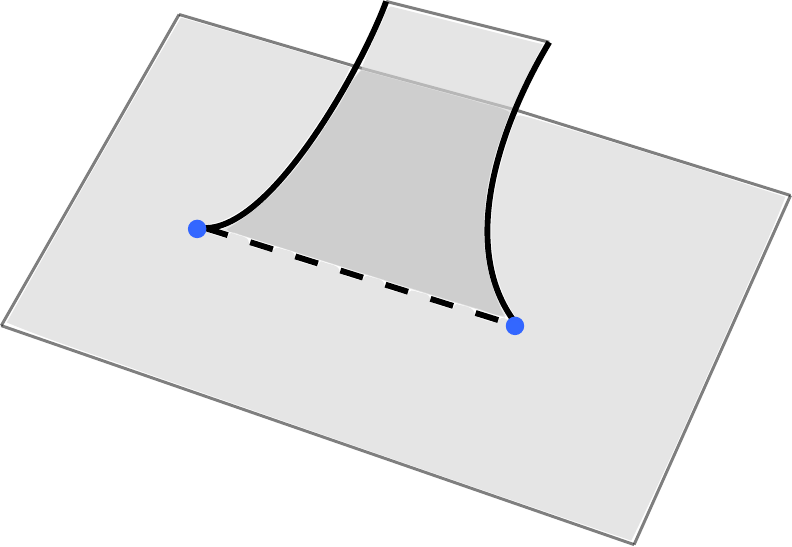}
    \caption{A flap.}
    \label{fig:flap}
    \end{subfigure}
    \begin{subfigure}[t]{0.3\linewidth}
    \includegraphics[scale=0.3]{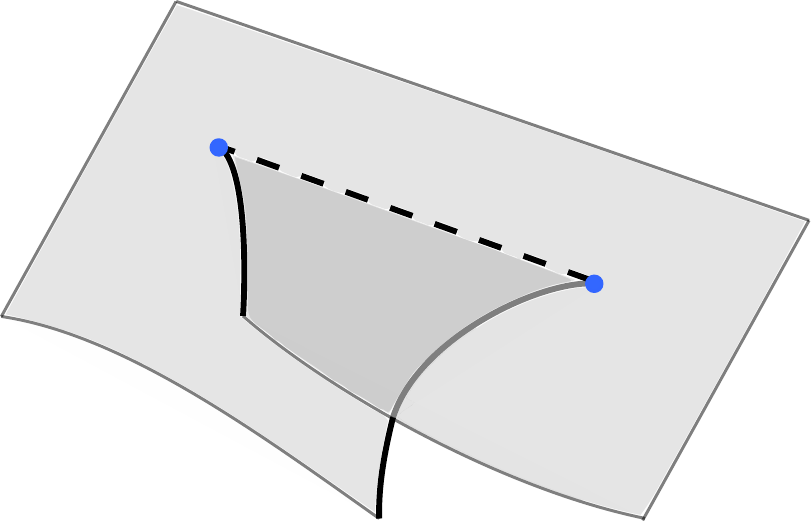}
    \caption{A pleat.}
    \label{fig:pleat}
    \end{subfigure}
    \caption{The figure shows (a) a local flap, (b) a flap, and (c) a pleat. The dashed line segments indicate the hyperbolic-regular values, whilst the thick black line segments indicate elliptic-regular values. The blue points indicate parabolic values.}
    \label{fig:flap_pleat}
\end{figure}

In the preimage under $\tau$ of a sufficiently small neighbourhood of a parabolic value, the bifurcation complex has two sheets: one sheet, the \emph{local base} $\mathcal B$, contains regular values and a compact line segment $\mathcal L$ of hyperbolic-regular values, and one sheet, the {\em local flap} $\mathcal F$, contains a line of elliptic-regular and of hyperbolic-regular values (which meet at a parabolic value) as well as regular values `between' these lines, see Figure \ref{fig:local_flap}. Both sheets intersect (or rather touch) each other along the line segment of hyperbolic-regular values including its parabolic end point. The topological boundary of $\mathcal F$ consists of the line segments of elliptic-regular and hyperbolic-regular values joint at the parabolic value and a line of regular values, called the {\em free boundary}.

Flaps and pleats/swallowtails now arise as follows: Consider a system with a compact line segment $\mathcal L$ of hyperbolic-regular values with parabolic end points denoted by $c_{1}$ and $c_2$. For $i \in \{ 1,2\}$, let $\mathcal{B}_{i}$ be their local bases and $\mathcal{F}_{i}$ their local flaps.
If one glues the free boundary of $\mathcal{F}_{1}$ to the free boundary of $\mathcal{F}_{2}$, this will define a \emph{flap topology} around $\mathcal L$, see Figure \ref{fig:flap}. If the free boundary of $\mathcal{F}_{1}$ is glued to the boundary of $\mathcal{B}_{2}$, and the free boundary of $\mathcal{F}_{2}$ is glued to the boundary of $\mathcal{B}_{1}$, this will define a \emph{pleat topology}, see Figure \ref{fig:pleat}. Efstathiou and Giacobbe \cite[Proposition 4]{Efstathiou2012} showed that the bifurcation complex in an open neighbourhood of $\mathcal L$ can have either the pleat topology or the flap topology. 

Efstathiou and Giacobbe \cite[Proposition 7]{Efstathiou2012} proved another interesting result: 
Let $p$ and $q$ be coprime integers and let $S^{3} := \{ (z_{1}, z_{2}) \in \C^{2} \mid \abs{z_{1}}^{2} + \abs{z_{2}}^{2} = 1 \}$ be the unit sphere in $\C^{2}$. Consider the (free) action of $\Z_{p} := \Z/p\Z$ on $S^{3}$ given by $(z_{1}, z_{2}) \mapsto \bigl(\exp(2 \pi i / p) z_{1}, \exp(2 \pi i q / p) z_{2}\bigr)$. The {\em lens space} $L(p,q) := S^{3} / \Z_{p}$ is the orbit space defined by this action. Then, with $\mathcal L$ as above, the type of lens space $L(p, 1)$ topologically embedded in $F^{-1}(\mathcal L)$ determines the monodromy of the Lagrangian fibration in a neighbourhood of $\mathcal L$ up to a sign determined by the choice of orientations.

\subsection{Subcritical Hamiltonian-Hopf bifurcations}\label{sec:subcritical_Hamiltonian-Hopf}

Recall from Section \ref{sec:supercritical_Hamiltonian-Hopf}, that a semitoric system with focus-focus points may arise via supercritical Hamiltonian-Hopf bifurcations from a toric one.
Analogously, a hypersemitoric system with flaps may arise from a semitoric one with focus-focus points via so-called \emph{subcritical Hamiltonian-Hopf bifurcations} by `replacing' a focus-focus point by a (small) flap, see for instance Dullin and Pelayo \cite{DullinPelayo2016}.

To be more precise, recall the normal form $\hat{H}_{s} = J+ X + s Y + a Y^{2}$ from Equation \eqref{eq:Hamiltonian-Hopf}: If the sign of $a$ is negative, then a focus-focus point (four complex eigenvalues) will first turn into a degenerate point (two purely imaginary eigenvalues of multiplicity $2$) and then will bifurcate into an elliptic-elliptic point (four purely imaginary eigenvalues) from the value of which, lying on a flap, two lines of elliptic-regular values emanate that connect the elliptic-elliptic value to the parabolic values (cf.\ Section \ref{sec:flaps_pleats}). The parabolic values are connected by a line of hyperbolic-regular values. 

In Figure \ref{fig:octagon_swallowtail}, an example of a semitoric system that went  through a subcritical Hamiltonian-Hopf bifurcation is displayed.

\subsection{Atoms, molecules, and classifications} \label{sec:hst_classification_molecules}

\begin{figure}[h]
    \centering
    \includegraphics[scale=0.55]{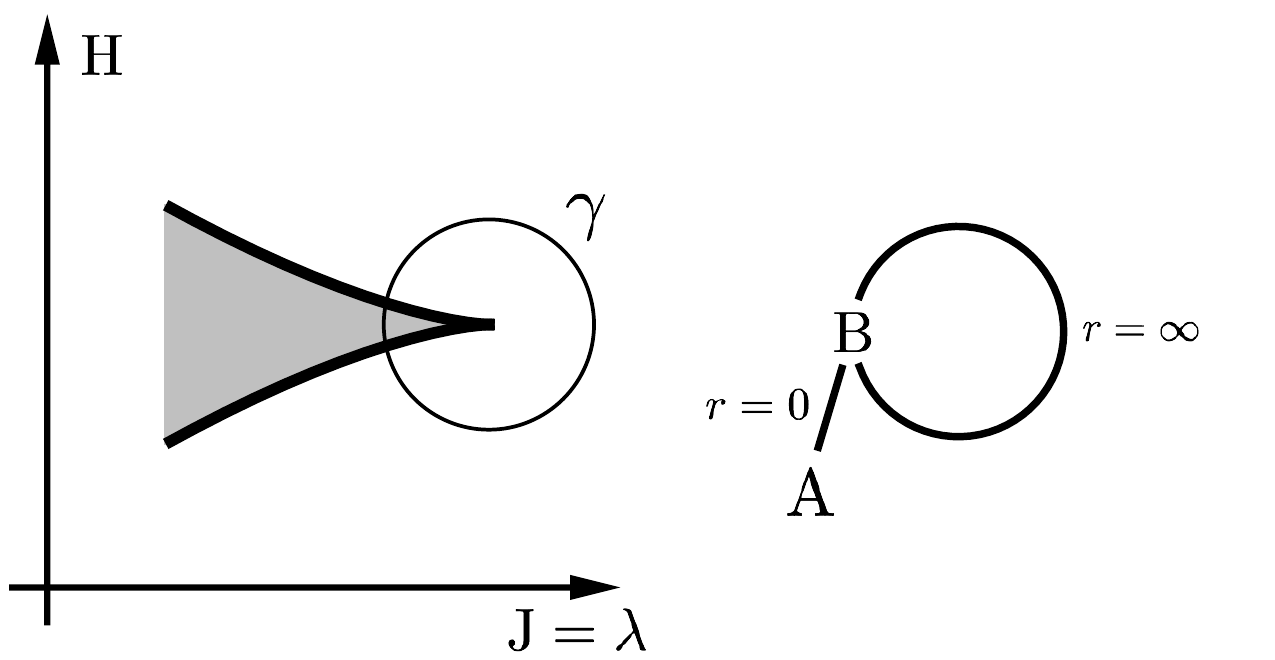}
    \caption{On the left, the bifurcation diagram of a parabolic singularity and a loop $\gamma$ around its critical value. On the right, the corresponding loop molecule.}
    \label{fig:loop_molecule_cusp}
\end{figure}

Recall from Section \ref{sec:AtomsMolecules} the notion of a marked molecule $W^{*}$, which is a complete isoenergy invariant of a 2 degree of freedom integrable system. The topology caused by the lines of elliptic-regular and hyperbolic-regular values in flaps and pleats (swallowtails) can in particular be described by marked molecules.
Here one can consider `loop molecules' (see Figure~\ref{fig:loop_molecule_cusp}) around the parabolic values with $B$-atoms describing the bifurcation of one of the two lines emanating from the cusp and $A$-atoms the other bifurcation.

The important result in this context is that the loop molecule around the cusp is uniquely defined and moreover `knows' what happens in its vicinity, in the sense that the loop 
molecule completely determines the topology of the corresponding singular torus fibration. This result directly follows from the fact that a single parabolic orbit (more precisely, the associated
compact singular fiber, which has the form of a cuspidal torus) gives rise to only one singular torus fibration up to a fibrewise homeomorphism, see Efstathiou and Giacobbe \cite{Efstathiou2012}. We conjecture that more is true in fact and there is only one
such torus fibration up to fibrewise diffeomorphisms, cf. Kudryavtseva and Martynchuk \cite{Kudryavtseva2021b}.

A similar topological result is known for elliptic-elliptic, elliptic-hyperbolic and focus-focus singularities of integrable systems on $4$-manifolds, but not so for hyperbolic-hyperbolic singularities (having multiple hyperbolic-hyperbolic points
on a singular fiber) which are in general not determined by their loop molecules only, see for instance \cite{Bolsinov2004, Bolsinov2006}. Interestingly, in the smooth case, the fibrewise classification turns out to be different also in the case of focus-focus singularities 
(having multiple points on the same singular fibre), see Bolsinov and Izosimov \cite{Bolsinov2019}.

The fibres of hypersemitoric systems will be classified by means of a `labeled graph' in the forthcoming paper by Gullentops and Hohloch \cite{Gullentops2023} which extends the special case of hyperbolic-regular fibres studied in Gullentops' thesis \cite{GullentopsThesis}. 

\subsection{Examples}

Hypersemitoric systems were first defined in Hohloch and Palmer \cite[Section 3]{Hohloch2021} who gave several examples for this class of systems. There are more examples in the paper by Gullentops and Hohloch \cite{Gullentops2022} and Gullentops' thesis \cite{GullentopsThesis}. 

\subsubsection{Hypersemitoric coupled angular momenta system}

\begin{figure}[ht]
    \centering
    \begin{subfigure}[t]{0.3\linewidth}
    \includegraphics[scale=0.45]{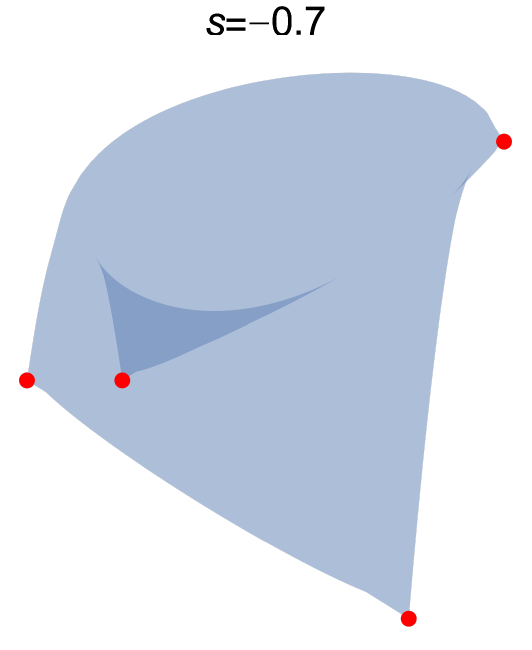}
    \caption{A flap.}
    \label{fig:hyp_cam_flap}
    \end{subfigure}
    \quad
    \begin{subfigure}[t]{0.3\linewidth}
    \includegraphics[scale=0.45]{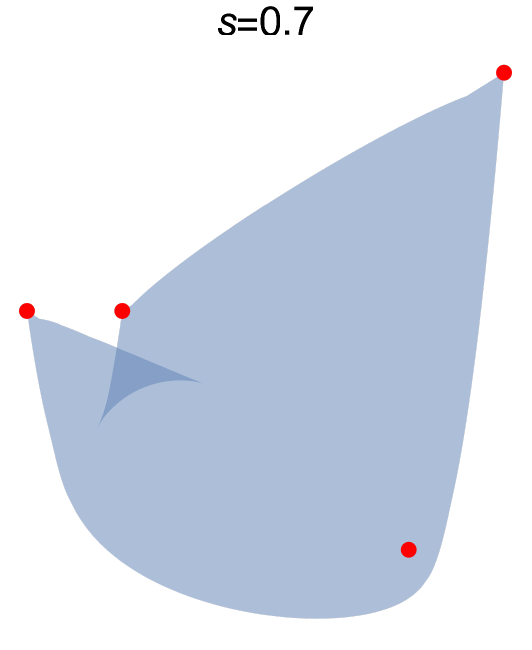}
    \caption{A pleat with a focus-focus point.}
    \label{fig:hyp_cam_pleat_ff}
    \end{subfigure}
    \quad
    \begin{subfigure}[t]{0.3\linewidth}
    \includegraphics[scale=0.45]{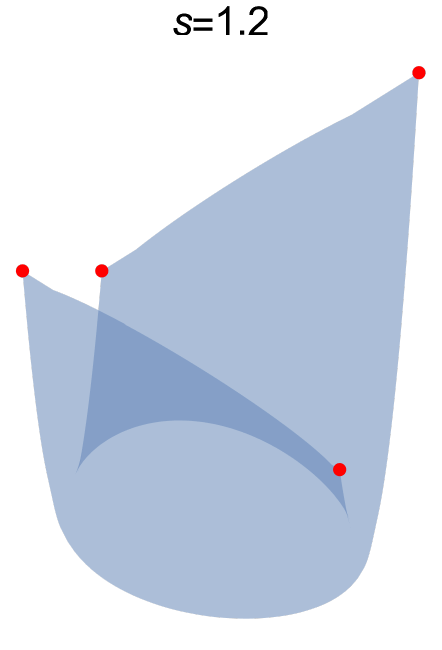}
    \caption{A pleat with an elliptic-elliptic point on it.}
    \label{fig:hyp_cam_pleat_ee}
    \end{subfigure}
    \caption{The momentum map image for the hypersemitoric coupled angular momenta system, for $R_{1} = 1$ and $R_{2} = 4$, and for selected values of $s$.}
    \label{fig:hyp_cam}
\end{figure}

Let $J$ and $H$ be as in the (semitoric) coupled angular momenta system, as discussed in Section \ref{sec:semitoric_cam}. We will now modify $H$, such that we instead consider the following:
\begin{equation*}
\tilde{H}(x_1,y_{1},z_{1},x_{2},y_{2},z_{2}) := H(x_1,y_{1},z_{1},x_{2},y_{2},z_{2}) + sz_{1}^{2},
\end{equation*}
with parameter $s \in \R$. Then, it turns out that the image of the momentum map $\tilde{F} = (J,\tilde{H})$, when the coupling parameter is $t = 0.5$ for which we always have a focus-focus value in the semitoric case (i.e.\ $s = 0$), we can generate flaps and pleats, see Figure \ref{fig:hyp_cam}. It turns out that the point $p_{1} = (0,0,1,0,0,-1)$ is of focus-focus type if $s_{p_{1}}^{-} < s < s_{p_{1}}^{+}$, where 
\begin{align*}
s_{p_{1}}^{\pm} = \frac{R_{1} \pm 2 \sqrt{R_{1} R_{2}}}{4R_{2}}.
\end{align*}
If $s < s_{p_{1}}^{-}$ or $s > s_{p_{1}}^{+}$, then $p_{1}$ is of elliptic-elliptic type. Numerics indicates that, if $R_{1} < R_{2}$, for $s < s_{p_{1}}^{-}$ a flap appears, and for some $s > s_{p_{1}}^{+}$, then a pleat appears. If $s \in \{s_{p_{1}}^{-},s_{p_{1}}^{+}\}$, then $(0,0,1,0,0,-1)$ is a degenerate singularity. This can be shown by a similar procedure as in Le Floch and Pelayo \cite[Proof of Proposition 2.5]{LeFlochPelayoCAM}. Furthermore, the point $p_2 = (0,0,-1,0,0,1)$ is a focus-focus point if $s_{p_{2}}^{-} < s < s_{p_{2}}^{+}$, where 
\begin{align*}
s_{p_{2}}^{\pm} = \frac{R_{1} \pm 2 \sqrt{R_{1} R_{2}} + 2R_{2}}{4R_{2}}.
\end{align*}
When $s < s_{p_{2}}^{-}$, then $\tilde{F}(p_{2})$ is an elliptic-elliptic value on the boundary of the momentum map image. For some $s > s_{p_{2}}^{+}$ we have that $\tilde{F}(p_{2})$ is an elliptic-elliptic value which joins the pleat created by $p_{1}$, see Figure \ref{fig:hyp_cam_pleat_ee}.

\subsubsection{The hypersemitoric octagon system}

\begin{figure}[h]
    \centering
    \includegraphics[scale=0.5]{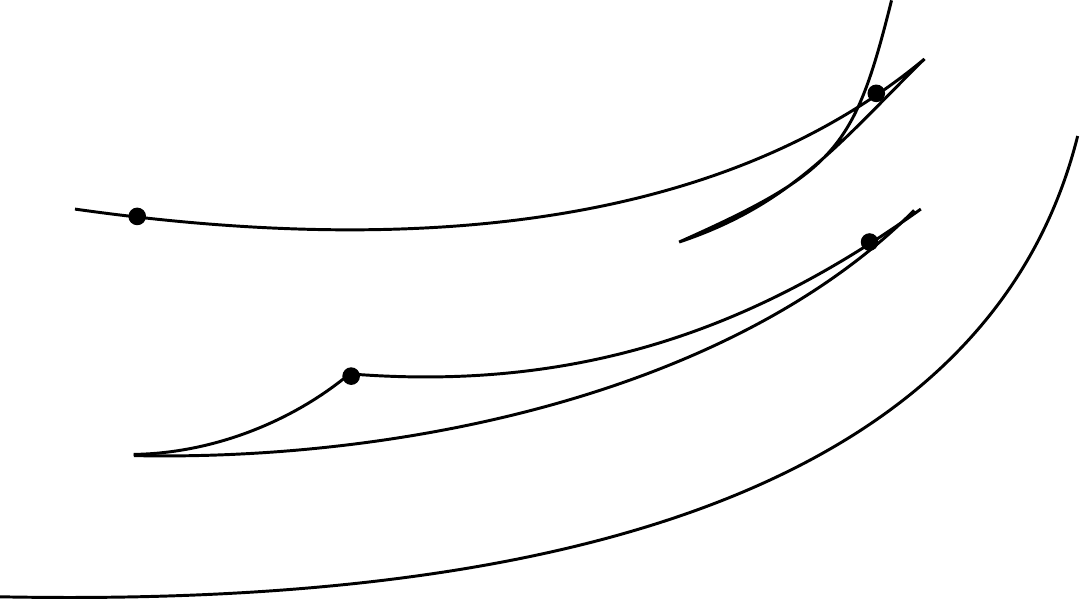}
    \caption{A sketch of a swallowtail and a flap appearing in the momentum map image of the hypersemitoric octagon system. The black points depict the rank $0$ singularities. For more details and plots, see Gullentops and Hohloch \cite{Gullentops2022}.}
    \label{fig:octagon_swallowtail}
\end{figure}

A specific family of examples can be created by taking the toric octagon system constructed in De Meulenaere and Hohloch \cite{Meulenaere2019} and, instead of perturbing it only to a semitoric system (cf.\ Section \ref{sec:semitoric_octagon}), more perturbation terms can be added to obtain a family of hypersemitoric systems. To be more precise, let $F=(J, H)$ be as in Section \ref{sec:semitoric_octagon} and modify $H$ to $H_{t}$ with $t = (t_{1}, t_{2}, t_{3}, t_{4}) \in \R^{4}$ via setting
\begin{equation*}
    H_{t} := (1 - 2t_{1})H + \sum_{i=1}^{4} t_{i} \gamma_{i},
\end{equation*}
with 
\begin{align*}
    &\gamma_{1}([z]) := \frac{1}{50} \Re\left( \bar{z}_{2}\bar{z}_{3}\bar{z}_{4}z_{6}z_{7}z_{8} \right), 
    &&\gamma_{2}([z]) := \frac{1}{50} \abs{z_{5}}^{4} \abs{z_{4}}^{4}, \\
    &\gamma_{3}([z]) := \frac{1}{50} \abs{z_{4}}^{4} \abs{z_{7}}^{4}, 
    &&\gamma_{3}([z]) := \frac{1}{50} \abs{z_{5}}^{4} \abs{z_{7}}^{4}.
\end{align*}
Gullentops and Hohloch \cite{Gullentops2022} proved the appearance of flaps and pleats/swallowtails and their collisions for certain values of the parameter $t$, see for example Figure \ref{fig:octagon_swallowtail}. Moreover, they studied the shape and topology for hyperbolic-regular fibres in the system $(J, H_t)$ and showed that, for fibres over a hyperbolic-regular value, not only double tori (`two tori stacked on top of each other' resp.\ a figure eight loop times $S^1$) are possible, but that the number of `tori stacked on top of each other' possibly appearing as fibre of a hyperbolic-regular value is bounded from above by 13.

\printbibliography

\end{document}